\documentclass[reqno]{amsart}
\usepackage[utf8]{inputenc}

\usepackage[utf8]{inputenc}

\usepackage{sky}
\usepackage{graphicx}
\usepackage{color}
\usepackage{caption}
\usepackage{halloweenmath}

\usepackage{tikz}
\usetikzlibrary{arrows.meta}
\usetikzlibrary{calc}

\newcommand{\Wg}{\mrm{Wg}}

\newcommand{\pr}{\mathbb{P}}

\newcommand{\Addresses}{{% additional braces for segregating \footnotesize
  \bigskip
  \footnotesize
  Ron Nissim, \textsc{Department of Mathematics, Massachusetts Institute of Technology, Cambridge, MA, USA.}\par\nopagebreak
  \textit{E-mail address}: \texttt{rnissim@mit.edu}

}}

\title{Sharp Estimates for Large N Weingarten Functions}
\author{Ron Nissim}
\date{}

\begin{document}

\begin{abstract}
    Weingarten functions provide a tool for computing Haar measure matrix integrals of polynomials in the matrix entries. An important property of Weingarten functions is their particularly simple large--$N$ limits. In 2017 Benoît Collins and Sho Matsumoto studied when this limit holds for Weingarten functions associated to integrals of products of $2n$ matrix entries, as $n \to \infty$, together with the matrix size $N$ \cite{Collins2017}. They showed that the large--$N$ limit is uniformly achieved as long as $n=o(N^{4/7})$, a result which already has applications to strong asymptotic freeness \cite{bordenave2023norm,bordenave2024strong}. However, their result is not optimal. They conjectured that their result should actually hold up to $n=o(N^{2/3})$ which is optimal. We prove this conjecture for the matrix groups $G \in \{\mathrm{U}(N)$, $\mathrm{O}(N)$, $\mathrm{Sp}(N)\}$. The proof proceeds by introducing a Markov process on permutations (pairings) which we call the unitary (orthogonal) \textit{Weingarten process}. We believe this process may have further applications to the theory of Weingarten functions. We also prove two new bounds regarding the large--$N$ limit of the Weingarten function in the regimes when $n=o(N^{4/5})$, and $n=o(N)$.
\end{abstract}
\maketitle

\section{Introduction and Main Results}

\subsection{Introduction}

For a compact matrix group  $G \in \{\mathrm{U}(N)$, $\mathrm{O}(N)$, $\mathrm{Sp}(N)\}$, let $\mu_G$ denote the  Haar probability measure on $G$, and let $U_{i,j}: G \to \C$ be the $i,j$ coordinate function. Then the integral $\int_G U_{i_1j_1} \cdots U_{i_m j_m}\overline{U}_{i'_1j'_1}\cdots \overline{U}_{i'_n j'_n} d\mu_G$ can be calculated in terms a class of special functions called Weingarten functions  \cite{weingarten1978asymptotic}. For general background on Weingarten functions we refer the reader to \cite{collins2022weingarten, collins2022moment,kostenberger2021weingarten}.

These Haar measure integrals which can be evaluated in terms of Weingarten functions have applications in many areas of probability and mathematical physics. These applications include strong convergence results in free probability \cite{collins2014strong,bordenave2023norm,collins2024free,bordenave2024strong}, analysis of famous unitary group integrals arising from mathematical physics such as the Harish-Chandra-Itzykson-Zuber integral and the Bars-Green/Brezin-Gross-Witten/Wadia integrals \cite{collins2003moments, novak2007truncations,collins2009asymptotics,goulden2011monotone1,goulden2011monotone2,novak2020complex,novak2022topological,collins2023tensor}, Yang-Mills lattice gauge theory \cite{novak20242d,CPS2023}, and quantum information theory \cite{collins2016random}. The theory of Weingarten functions has been particularly successful for studying large--$N$ limits of matrix models in these areas due to the particularly simple form of Weingarten functions in the large--$N$ limit. For instance, the unitary Weingarten function takes in a permutation input $\sigma \in S_n$, and we have the following $N$ limit formula: for $n$ fixed,
\begin{equs}\label{Large N limit of Weingarten Functions}
    \Wg_N^{\mathrm{U}}(\sigma) = N^{-n-|\sigma|} \mathrm{Moeb}(\sigma)(1+O(N^{-2})),
\end{equs}
where for $\sigma \in S_n$ with $l$ cycles of length $C_1,...,C_l$, $|\sigma|:=n-l$, and $\mathrm{Moeb}(\sigma):=(-1)^{|\sigma|}\prod_{i=1}^{l} \mathrm{Cat}(C_l-1)$, for $\mathrm{Cat}(m):= \frac{\binom{2m}{m}}{m+1}$ defined as the $m$-th Catalan number. $\mathrm{Moeb}(\sigma)$ is called the Moebius function of $\sigma$.

For many applications listed above, it is important to understand Weingarten functions in the case when both $n\to \infty$ and $N \to \infty$. In particular, it is natural to ask under what conditions does $\Wg_N^{\mathrm{U}}(\sigma) = N^{-n-|\sigma|} \mathrm{Moeb}(\sigma)(1+o(1))$ hold uniformly for $\sigma \in S_n$ as $n \to \infty$ and $N \to \infty$ simultaneously. This question is more delicate than the original large--$N$ limit question answered by \eqref{Large N limit of Weingarten Functions}, however, there have been results in this direction. In particular, Collins and Matsumoto \cite{Collins2017}, proved that for any $N \geq \sqrt{6}n^{7/4}$, we have
\begin{equs}\label{Collins Matsumoto bound}
    \frac{1}{1-\frac{n-1}{N^2}} \leq \frac{N^{-n-|\sigma|} \Wg_N^{\mathrm{U}}(\sigma)}{\mathrm{Moeb}(\sigma)} \leq \frac{1}{1-\frac{6n^{7/2}}{N^2}},
\end{equs}
for all $\sigma \in S_n$. 

However, Collins and Matsumoto remark that the exponent in the upper bound of \eqref{Collins Matsumoto bound} is not tight. For the full cycle $\sigma=(1 \, 2 \, \dots \, n) \in S_n$, $\Wg_N(\sigma)$ has an explicit formula which can be analyzed to deduce that
\begin{equs}
    \frac{N^{-n-|\sigma|} \Wg_N^{\mathrm{U}}(\sigma)}{\mathrm{Moeb}(\sigma)} \to 1
\end{equs}
 as long as $n,N \to \infty$ with $n = o(N^{2/3})$. This suggests a uniform bound of the form 

\begin{equs}\label{Our Main Bound Intro}
    \frac{N^{-n-|\sigma|} \Wg_N^{\mathrm{U}}(\sigma)}{\mathrm{Moeb}(\sigma)} \leq \frac{1}{1-\frac{Cn^{3}}{N^2}},
\end{equs}
for some constant $C>0$ and $n < C^{-1}N^{2/3}$. This is indeed the statement of our main theorem, Theorem \ref{Main Weingarten Bound}.  We also prove the analogous result for the $\mathrm{O}(N)$ and $\mathrm{SP}(N)$ case, Theorem \ref{Main Bound Orthogonal}.

The starting point to proving \eqref{Our Main Bound Intro} is the path expansion formula for the Weingarten function proven in \cite{Collins2017},
\begin{equs}\label{Path Expansion Intro}
    (-1)^{|\sigma|} N^{n+|\sigma|}\mathrm{Wg}^{\mathrm{U}}_N(\sigma) = \sum_{g=0}^{\infty} |\mathrm{P}(\sigma, |\sigma|+2g)| N^{-2g},
\end{equs}
where $\mathrm{P}(\sigma, |\sigma|+2g)$ denotes the set of all directed paths of length $n+|\sigma|+2g$ from $\sigma \to \emptyset$ in the Weingarten graph which will be defined in Section \ref{U Weingarten Graph Subsection}. Then the crux of the argument will be to introduce a Markov process which generates uniformly random paths in $\mathrm{P}(\sigma, |\sigma|)$ which we call the \textit{Weingarten Process}. The crux of the proof will then reduce to bounding certain statistics of this process. This process is a weighted version of the uniform transposition shuffle which has been extensively studied studied \cite{diaconis1981generating,schramm2005compositions,berestycki2006phase,jain2024hitting}. The connection between matrix integrals and transposition shuffles has appeared in previous works including \cite{levy2008schur,park2023wilson}.

In Section \ref{Small Permutations Section} of the paper we obtain bounds on the Weingarten function for a wider set of $(n,N)$ values. Theorem \ref{n^{5/4} Theorem} is a logarithmic version of the large--$N$ limit result which holds up to the $n=o(N^{4/5})$ regime. However, Theorem \ref{n^{5/4} Theorem}, does not provide a uniform error over all permutations, and is strongest for permutations which have at least one large cycle. The proof of Theorem \ref{n^{5/4} Theorem} follows from analyzing the loop equation satisfied by the Weingarten function. Lastly, Theorem \ref{Main Weingarten Bound} and Theorem \ref{n^{5/4} Theorem} both provide bounds which are made to handle permutations $\sigma$, with $|\sigma|$ relatively large. We also study the ``small permutation'' regime in Section \ref{Small Permutations Section}  when $|\sigma| \leq k$ for some constant $k$ and $n, N \to \infty$. For such permutations, we can prove a bound for an even larger range of $(n,N)$ values.  More precisely, if $n,N \to \infty$, with $n = o(N)$ we  prove that $N^{-n-|\sigma|} \Wg_N^{\mathrm{U}}(\sigma) \to \mathrm{Moeb}(\sigma)$. The more precise bound is stated in Theorem \ref{Small Permutations Theorem}. The proof of this result is elementary and purely combinatorial, again using the path expansion \eqref{Path Expansion Intro}. 

%Unfortunately, proving any sort of bounds in the so called \textit{unstable} regime when $n>N$, almost certainly requires a completely different approach as the path expansion as \eqref{Path Expansion Intro} no longer converges in this regime. There are still some explicit formulas in the unstable regime. For instance, the Weingarten function $\Wg_N^{\mathrm{U}(N)}(\sigma)$ has the representation theoretic formula
%\begin{equs}
%    \Wg_N^{\mathrm{U}(N)}(\sigma) = \frac{1}{n!} \sum_{\substack{\lambda \vdash n, \\ \ell(\lambda) \leq N}} [\chi_{\lambda}(\groupid)\chi_{\lambda}(\sigma)\prod_{(i,j) \in \lambda}(N+j-i)^{-1}],
%\end{equs}
%where $\chi_{\lambda}$ are the characters of the irreducible representations of $S_n$ indexed by partitions $\lambda \vdash n$, and $\ell(\lambda)$ denotes the number of parts in $\lambda$ \cite{collins2006integration}. We leave some open questions,...

\vspace{3mm}
\textbf{Acknowledgements: } The author would like to thank Sky Cao and Jacopo Borga for reading over an earlier draft and providing useful comments and suggestions. The author would also like to thank the reviewer for many helpful comments and suggestions. Lastly, the author would like to thank Jacopo Borga, Sky Cao, and Scott Sheffield for many helpful and inspiring discussions. The author is supported by the NSF GRFP-2141064.

\subsection{Main Results and Discussion}

We encourage the reader to refer to Subsection \ref{U Weingarten Function Section} for notation used in this section.

The main result of the paper is the following tight uniform bound on the Unitary Weingarten function conjectured by Collins and Matsumoto \cite{Collins2017}.

\begin{theorem}\label{Main Weingarten Bound}
    Fix $C=6\sqrt{8}\times10^6$. Then for $N \geq \sqrt{C}n^{3/2}$ and all $\sigma \in S_n$ we have,
    \begin{equs}
        \frac{N^{n+|\sigma|}\Wg_N^{\mathrm{U}}(\sigma)}{\mathrm{Moeb}(\sigma)} \leq \frac{1}{1-\frac{Cn^3}{N^2}}.
    \end{equs}
\end{theorem}

\begin{remark}
     For the full cycle $\sigma_n:= (1 \, 2 \, \dots \, n)$, 
     \begin{equs}\label{Full Cycle Weingarten Formula}
         \Wg_N^{\mathrm{U}}(\sigma_n)=\frac{(-1)^{n-1}\mathrm{Cat}(n-1)}{(N-n+1)\dots(N-1)N(N+1)\dots(N+n-1)}
     \end{equs}
     \cite[Proposition 6.1]{collins2009some}.
     As a result, it can be checked that the bound from Theorem \ref{Main Weingarten Bound} is tight in the sense that the normalized Weingarten function is no longer approximated by the Moebius function once $n$ is comparable to $N^{2/3}$.
\end{remark}
The proof of Theorem \ref{Main Weingarten Bound} will be given in Section \ref{Proof of Main theorem section}.

We also prove a similar bound for $\mathrm{O}(N)$ and $\mathrm{SP}(N)$ in Section \ref{Orthogonal Section}. In these cases the Weingarten function will be defined in Section \ref{Orthogonal Section}, Lemma \ref{Orthogonal Weingarten Expansion} and Remark \ref{Orthogonal to Symplectic}, and takes in pairings $\pi \in \mathcal{P}_2(2n)$. The bound for $\mathrm{O}(N)$ and $\mathrm{SP}(N)$ is stated below. 

\begin{theorem}\label{Main Bound Orthogonal}
    For $G \in \{\mathrm{O}(N),\mathrm{SP}(N)\}$, and $N \geq 2 \times 10^6 n^{3/2}$, we have
    \begin{equs}
        \frac{N^{-n-|\pi|}\mathrm{Wg}_N^{G}(\pi)}{\mathrm{Moeb}(\pi)}\leq \frac{1}{(1-10^6n^{3/2}/N)^2},
    \end{equs}
    for any $\pi \in \mathcal{P}_2(2n)$.
\end{theorem}

Moving back to the $G=\mrm{U}(N)$ setting, in Section \ref{Small Permutations Section} we will establish two quantitative large $N$ results that will work for an even larger regime of $(n,N)$ values than Theorem \ref{Main Weingarten Bound}, but will no longer be uniform over permutations $\sigma \in S_n$. The following result addresses the case when $n=o(N)$, but $|\sigma|$ remains bounded as $n,N \to \infty$.
\begin{theorem}\label{Small Permutations Theorem}
    For $N > \sqrt{48}en$, and some constant $C_{|\sigma|}$ depending only on $|\sigma|$,
    \begin{equs}
        |N^{n+|\sigma|}|\Wg_N^{\mathrm{U}(N)}(\sigma)-\mathrm{Moeb}(\sigma)| \leq \frac{C_{|\sigma|}n^2}{N^2-48en^2}.
    \end{equs}
    Where $e$ refers to Euler's constant $e=2.718\dots$
\end{theorem}

\begin{remark}
    From the proof of Theorem \ref{Small Permutations Theorem}, one can choose $C_{|\sigma|} \leq C e^{|\sigma|^2}$ for some constant $C$ independent of $|\sigma|$, and if one is more careful, in fact $C_{|\sigma|} \leq C |\sigma|^{|\sigma|}$ should also work.
\end{remark}

We also obtain yet another result which studies the logarithmic version of the large--$N$ limit motivated as follows. From the formula for the full cycle \eqref{Full Cycle Weingarten Formula}, we still have the correct limiting asymptotics in a logarithmic sense as long as $n=o(N)$, that is,
\begin{equs}
    \log\bigg( \frac{N^{n+|\sigma|}|\Wg_N^{\mathrm{U}(N)}(\sigma)|}{|\mathrm{Moeb}(\sigma)|}\bigg) \to 1
\end{equs}
as $N,n \to \infty$ with $n=o(N)$. 

This leads us to the following theorem which deals with the $n=o(N^{5/4})$ regime.

\begin{theorem}\label{n^{5/4} Theorem}
    Suppose $n \leq 10^{-4}N^{4/5} $. Then for $\sigma \in S_n$ with $|\sigma|\geq 4$,
    \begin{equs}       
    \log\bigg(\frac{N^{n+|\sigma|}|\Wg_N^{\mathrm{U}(N)}(\sigma)|}{|\mathrm{Moeb}(\sigma)|}\bigg)\leq 25\frac{n^{5/2}}{N^2}|\sigma|+\log|\sigma|+2.
    \end{equs}  
\end{theorem}

\begin{remark}
    Recall that $\mathrm{Cat}(n) \sim \frac{4^{n}}{\sqrt{\pi}n^{3/2}}$, so as long as $N \gg n^{5/4}$, and $\sigma$ is permutation with at least one large cycle relative to $|\sigma|$, Theorem \ref{n^{5/4} Theorem} can be interpreted as, 
    \begin{equs}
        \log(N^{n+|\sigma|}|\Wg_N^{\mathrm{U}(N)}(\sigma)) = \log|\mathrm{Moeb}(\sigma)|+\text{Lower Order Terms}.
    \end{equs}
\end{remark}

Checking against the formula for the Weingarten function evaluated at the full cycle as in the prior discussion, Theorem \ref{n^{5/4} Theorem} seems to be sub-optimal. One would expect a similar result to hold in the entire $n=o(N)$ regime. We leave this as a problem for future work.

\section{Notation and Background}

\subsection{Partitions and Permutations}\label{PPP Subsection}

A partition of $n$, denoted $\lambda \vdash n$, is a sequence of numbers $\lambda:=(\lambda_1,\lambda_2,\dots,\lambda_l)$ with $\lambda_1 \geq \lambda_2 \geq \dots \geq \lambda_l \geq 1$, and $\lambda_1+\lambda_2+\dots +\lambda_l=n$. For a partition $\lambda \vdash n$ with $\lambda=(\lambda_1,\lambda_2,\dots,\lambda_l)$, we define $\ell(\lambda)=l$ to be the length of the partition. The letter $l$ will also frequently be used to denote the number of cycles in a permutation since the cycle structure of the permutation will then be a partition $\lambda$ with $\ell(\lambda)=l$.

Let $[n]:=\{1,2,\dots,n\}$. We let $S_n$ denote the group of permutations acting on $[n]$ (the set of bijective functions $\sigma:[n] \to [n]$). A permutation $\sigma \in S_n$ is a \textit{cycle} if there is no nontrivial subset $X \subset [n]$ such that $\sigma(X)=X$. A permutation $\tau \in S_n$ is a \textit{transposition} if it fixes all but two elements of $[n]$. Recall the following facts about permutations.

\noindent \textbf{Fact 1.} Any permutation $\sigma \in S_n$ can be factored uniquely up to ordering into a product of disjoint cycles. The lengths of these cycles $\lambda_1 \geq \lambda_2,\dots \geq\lambda_l$, uniquely determine a partition $\lambda \vdash n$ called the cycle type of $\sigma$. We let $S_{\lambda} \subset S_n$ denote the collection of all permutations with cycle type $\lambda$.

\noindent \textbf{Fact 2.} For any permutation $\sigma \in S_{\lambda}$, the conjugacy class of $\sigma$, $K_{\sigma}:=\{\Tilde{\sigma}=\eta^{-1}\sigma \eta: \eta \in S_n\}$, is precisely equal to $S_{\lambda}$.

\noindent \textbf{Fact 3.} Any permutation $\sigma \in S_n$ can be factored into a product of transpositions $\sigma=\tau_1 \dots \tau_m$. The minimal number of transpositions necessary to factor $\sigma$ this way is denoted by $|\sigma|$, and letting $l$ denote the number of cycles in the disjoint cycle decomposition of $\sigma$, $|\sigma|=n-l$.

\noindent For any $\sigma \in S_n$, and a permutation $\tau = (i \, j)$, 
\begin{equs}
    |\tau \sigma|= 
    \begin{cases}
        |\sigma|-1 \text{   if } i \text{ and } j \text{ belong to the same cycle of } \sigma.\\
        |\sigma|+1 \text{   otherwise.}
    \end{cases}
\end{equs}

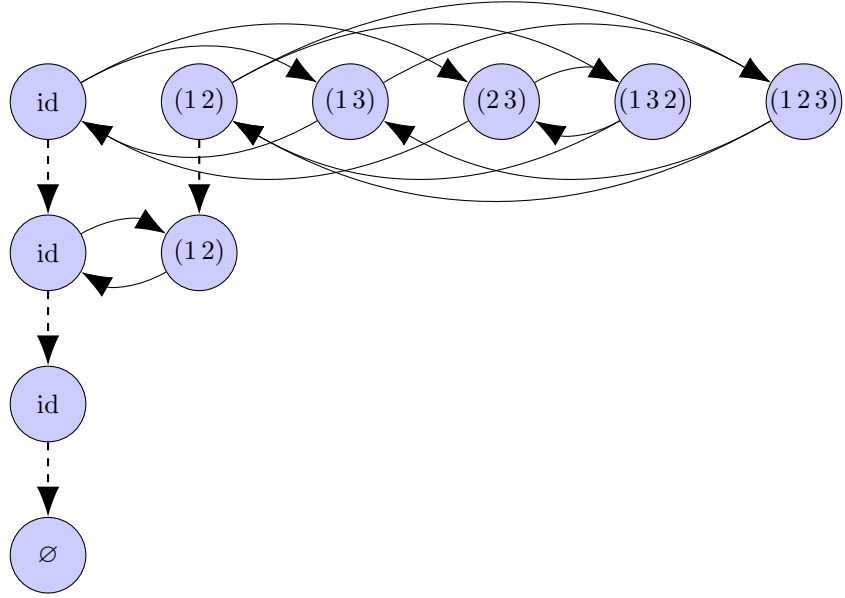
\begin{figure}
    \centering
    \begin{tikzpicture}[
    ->, % Arrow style for directed edges
    >={Latex[length=4mm, width=3mm]}, % Larger arrowheads
    node/.style={circle, draw, fill=blue!20, minimum size=10mm, inner sep=0pt, align=center},
    every edge/.style={draw, thick}
] % Base style for edges
    % Nodes
    \node[node](p1) at (0,6) {$\groupid$};
    \node[node](p2) at (2,6) {$(1\, 2)$};
    \node[node](p3) at (4,6) {$(1\, 3)$};
    \node[node](p4) at (6,6) {$(2\, 3)$};
    \node[node](p5) at (8,6) {$(1\, 3 \, 2)$};
    \node[node](p6) at (10,6) {$(1\, 2\,  3)$};
    \node[node](pp1) at (0,4) {$\groupid$};
    \node[node](pp2) at (2,4) {$(1\, 2)$};
    \node[node](ppp1) at (0,2) {$\groupid$};
    \node[node](pppp1) at (0,0) {$\varnothing$};
    \draw (p1) edge[dashed] (pp1);
    \draw (p2) edge[dashed] (pp2);
    \draw (pp1) edge[dashed] (ppp1);
    \draw (ppp1) edge[dashed] (pppp1);
    \draw (p1) to[bend left=30] (p3);
    \draw (p1) to[bend left=30] (p4);
    \draw (p3) to[bend left=30] (p1);
    \draw (p4) to[bend left=30] (p1);
    \draw (p2) to[bend left=30] (p5);
    \draw (p5) to[bend left=30] (p2);
    \draw (p2) to[bend left=30] (p6);
    \draw (p6) to[bend left=30] (p2);
    \draw (p3) to[bend left=30] (p6);
    \draw (p6) to[bend left=30] (p3);
    \draw (p4) to[bend left=30] (p5);
    \draw (p5) to[bend left=30] (p4);
    \draw (pp1) to[bend left=30] (pp2);
    \draw (pp2) to[bend left=30] (pp1);
\end{tikzpicture}

\caption{The figure above displays the unitary Weingarten graph $\mathcal{W}^{\mathrm{U}}$ restricted to layers 0 through 3, $\bigsqcup_{n=0}^{3} S_n$.}
\label{fig:Weingarten Graph}
\end{figure}

\subsection{The Unitary Weingarten Graph}\label{U Weingarten Graph Subsection}

The unitary Weingarten graph $\mathcal{W}^{\mathrm{U}}$, introduced in \cite{Collins2017}, has vertex set $\bigsqcup_{n \geq 0} S_n$ ($S_0:= \{\varnothing\}$). Each copy of $S_n$ is referred to as the ``$n$-th level'' of the graph. In each level, there is a directed solid edge going from $\sigma \to \sigma'$ if there is a transposition of the form $(i \, n)$ such that $\sigma' = (i \, n) \sigma$. If a permutation $\sigma \in S_n$ fixes $n$, then there is a dashed line connecting $\sigma \to \sigma|_{[n-1]} \in S_{n-1}$. See Figure \ref{fig:Weingarten Graph} for a picture of the unitary Weingarten graph. Lastly, $\mathrm{P}(\sigma, l)$ will denote the set of paths from $\sigma \to \varnothing$ which take exactly $l$ solid directed edges. Note this number is only non-zero if $l \geq |\sigma|$. 
 
\subsection{The Unitary Weingarten Function}\label{U Weingarten Function Section}

Collins and Matsumoto proved the following large $N$ expansion for the unitary Weingarten function, which we will simply take as our definition since it is the starting point of our analysis.

\begin{definition}\cite[Lemma 2.5]{Collins2017}\label{Weingarten Definition}
    For $n \leq N$, the unitary Weingarten function $\Wg_N^{\mathrm{U}(N)}: S_n \to  \C$ can be defined as follows,
    \begin{equs}\label{Path Expansion Unitary}
        \mathrm{Wg}^{\mathrm{U}}_N(\sigma) = (-1)^{|\sigma|} N^{-n-|\sigma|}\sum_{g=0}^{\infty} |\mathrm{P}(\sigma, |\sigma|+2g)| N^{-2g}.
    \end{equs}
    Moreover we define the normalized Weingarten function,
    \begin{equs}
        \overline{\mathrm{Wg}}^{\mathrm{U}}_N(\sigma):=N^{n+|\sigma|}\mathrm{Wg}^{\mathrm{U}}_N(\sigma)
    \end{equs}
\end{definition}
Next we recall the definition of the Moebius function of a permutation. 
\begin{definition}\label{Moebius Def}
    Suppose $\sigma \in S_n$ has cycle decomposition $\sigma=\sigma_1\cdots \sigma_l$, where each $\sigma_i$ is a cycle of length $C_i$, then the \textit{Moebius function} of $\sigma$ is,
    \begin{equs}
        \mathrm{Moeb}(\sigma):= (-1)^{|\sigma|}\prod_{i=1}^{l}\mathrm{Cat}(C_i-1),
    \end{equs}
    where $\mathrm{Cat}(k):= \frac{(2k)!}{k!(k+1)!}$ is the $k$th Catalan number.
\end{definition}

The large--$N$ limit of the normalized Weingarten function is the Moebius function in the following sense.

\begin{lemma}\cite[Lemma 3.3]{Collins2017}\label{Path Enumeration Formula}
    We have that,
    \begin{equs}
        \mathrm{P}(\sigma,|\sigma|) = |\mathrm{Moeb}(\sigma)|,
    \end{equs}
    so as a result for any fixed $n$ and $\sigma \in S_n$,
    \begin{equs}
        \lim_{N \to \infty} \overline{\mathrm{Wg}}^{\mathrm{U}}_N(\sigma)=\mathrm{Moeb}(\sigma).
    \end{equs}
\end{lemma}

The main application of the unitary Weingarten function is the following unitary Haar measure integration formula,
\begin{theorem}\cite[Theorem 4.4]{collins2022weingarten}
    For any indices $\mathbf{i}=(i_1,\dots,i_n)$, $\mathbf{j}=(j_1,\dots,j_n)$, $\mathbf{i}'=(i_1',\dots,i_n')$, and $\mathbf{j}'=(j_1',\dots,j_n')$, we have,
    \begin{equs}
        \int_{\mathrm{U}(N)} U_{i_1,j_1} \cdots U_{i_n,j_n} \overline{U}_{i_1',j_1'} \cdots \overline{U}_{i_n',j_n'} dU = \sum_{\sigma_1, \sigma_2 \in S_n} \Wg_N^{\mathrm{U}(N)}(\sigma_1^{-1}\sigma_2) \delta_{\sigma_1(\mathbf{i})=\mathbf{i}'}\delta_{\sigma_2(\mathbf{j})=\mathbf{j}'}.
    \end{equs}
\end{theorem}
\begin{remark}
    The large $N$ expansion can actually be derived from the integration relation above, via iteration of a loop equation for unitary integration \cite{Collins2017}.
\end{remark}

For further background on the unitary Weingarten function, we refer the reader to \cite{collins2022weingarten, collins2022moment,kostenberger2021weingarten}. 

\section{Unitary Group: Proof of Theorem \ref{Main Weingarten Bound}}\label{Proof of Main theorem section}

The reader is encouraged to briefly review Subsections \ref{PPP Subsection} and \ref{U Weingarten Graph Subsection}, and Definition \ref{Weingarten Definition} before reading this section.

In order to prove Theorem \ref{Main Weingarten Bound}, we will need to introduce a Markovian way to sample paths uniformly from $\mathrm{P}(\sigma,|\sigma|)$ and $\bigcup_{\sigma \in S_{\lambda}}\mathrm{P}(\sigma,|\sigma|)$, a class of paths appearing in the path expansion formula for the Weingarten function \eqref{Path Expansion Intro}. This is the content of Section \ref{Unitary Weingarten Process Section}.

\subsection{The Unitary Weingarten Process}\label{Unitary Weingarten Process Section}

The following enumeration formula for $|P(\sigma,|\sigma|)|$ will be crucial for us to define the Weingarten Process.

Before defining the Weingarten process, recall from subsection \ref{PPP Subsection}, for a partition $\lambda \vdash n$, $S_{\lambda}:= \{\sigma \in S_n: \sigma \text{ has cycle structure } \lambda\}$. If $\sigma \in S_{\lambda}$, we will denote $\lambda:=\mathrm{shape}(\sigma)$.

\begin{definition}
    The unitary Weingarten process with initial data $\sigma \in S_{n}$, $\mathrm{WP}(\sigma)$ is a random process on $\mathrm{P}(\sigma, |\sigma|)$ defined by the following algorithm.
    \begin{enumerate}
        \item Let $\sigma_0 =\sigma$, and $n_0=n$. Set $k=0$.
        \item If: $\sigma_k$ does not fix $n_k$, move to the next step.\\
        Else if: $n_k>1$ and $\sigma_k$ fixes $n_k$, set $n_{k+1}=n_k-1$, $\sigma_{k+1}=\sigma_k|_{[n_{k+1}]} \in S_{n_{k+1}}$, and repeat this step with $k$ now set to $k+1$. \\
        Else if: $n_k=1$, terminate the process.
        \item Suppose the cycle in the disjoint cycle decomposition of $\sigma_k$ containing $n_k$ is $(i_1\,...\, i_{\ell} \, n_k)$. With probability $\frac{\mathrm{Cat(j-1)\mathrm{Cat}(\ell-j)}}{\mathrm{Cat}(\ell)}$, let $\sigma_{k+1}=\sigma_k (i_j \, n_k) \in S_{n_k}$, $n_{k+1}=n_k$, and now set $k$ to $k+1$ and go back to step 2.
    \end{enumerate}
\end{definition}

\begin{definition}
    For a partition $\lambda \vdash n$, let $\mathrm{WP}(\lambda)$ be a process generating a random path of permutations $(\sigma_0,\sigma_1,...,\sigma_{2n-\ell(\lambda)})$ as follows:
    \begin{enumerate}
        \item Take $\sigma \sim \mathrm{Uniform}(S_{\lambda})$.
        \item Generate a random path $(\sigma_0,..., \sigma_{2n-\ell(\lambda)}) \in \mathrm{P}(\sigma,|\sigma|)$ according to $\mathrm{WP}(\sigma)$.
    \end{enumerate}
    Defining $\lambda_i := \mathrm{shape}(\sigma_i)$, from $\mathrm{WP}(\lambda)$, we also obtain a random sequence of partitions, $(\lambda_0,\lambda_1,...,\lambda_{2n-\ell(\lambda)})$.
\end{definition}

\noindent $\mathrm{WP}(\sigma)$ and $\mathrm{WP}(\lambda)$ come with the canonical filtration $(\mathcal{F}_k)_{k \leq n+|\sigma|}$ where $\mathcal{F}_k$ is the sigma algebra generated by $\{\sigma_j\}_{j \leq k}$.

\begin{lemma}\label{Unitary Process is uniform over paths}
    For any permutation $\sigma \in S_{\lambda}$,
    \begin{equs}
        \mathrm{WP}(\sigma) \sim \mathrm{Uniform}(\mathrm{P}(\sigma, |\sigma|)),
    \end{equs}
    and for any partition $\lambda \vdash n$,
    \begin{equs}
        \mathrm{WP}(\lambda) \sim \mathrm{Uniform}\bigg( \bigcup_{\pi \in S_{\lambda}} \mathrm{P}(\pi, |\pi|)\bigg).
    \end{equs}
\end{lemma}

\begin{proof}
    The proof is via induction on $|\sigma|$. The base case $|\sigma|=0$ is obvious. 

    For the inductive step, suppose first that $\sigma$ does not fix $n$. Then let $(i_1 \, ...\, i_{\ell} \, n)$ be the cycle containing $n$, and the first step in $\mathrm{P}(\sigma,|\sigma|)$ must take $\sigma \to \sigma (i_j \, n)$ for some $i_j$. Thus
    \begin{equs}
        \mathrm{P}(\sigma,|\sigma|) \simeq \bigcup_{j=1}^{l}\mathrm{P}(\sigma (i_j \, n),|\sigma|-1).
    \end{equs}
    In other words, using the enumeration formula in Lemma \ref{Path Enumeration Formula}, when sampling $p \in \mathrm{P}(\sigma,|\sigma|)$ uniformly, with probability
    \begin{equs}
        \frac{|\mathrm{P}(\sigma (i_j \, n),|\sigma|-1)|}{|\mathrm{P}(\sigma,|\sigma|)|}=\frac{\mathrm{Cat(j-1)\mathrm{Cat}(\ell-j)}}{\mathrm{Cat}(\ell)},
    \end{equs}
    the first step of the path is $\sigma \to \sigma (i_j \, n)$, followed by a uniformly random path in $\mathrm{P}(\sigma (i_j \, n),|(i_j \, n)\sigma|)$. But the Weingarten process obeys this recursion by construction.
\end{proof}

\subsection{Markovian Properties of the Weingarten Process}

It is a direct consequence of the definition of the Weingarten process, that $(\mathrm{WP}(\sigma)$ is a homogeneous Markov process with respect to the filtration $(\mc{F}_{k})_{k \geq 1}$.

\begin{lemma}\label{Markov Property}
    $(\mathrm{WP}(\sigma)$ and $(\mathrm{WP}(\lambda)$ are homogeneous Markov processes with respect to the filtration $(\mc{F}_{t})_{t \geq 1}$. That is, for any $\pi \in \bigcup_{n=0}^{\infty} S_n$,
    \begin{equs}
        \pr(\sigma_{t+1}=\pi| \mathcal{F}_{t})=\pr(\sigma_{t+1}=\pi|\sigma_t),
    \end{equs}
    and for any permutations $\pi_1,\pi_0$,
    \begin{equs}
        \pr(\sigma_{t+1}=\pi_1| \sigma_{t}=\pi_0)= \pr(\sigma_{1}=\pi_1| \sigma_{0}=\pi_0).
    \end{equs}
\end{lemma}

Moreover, it is a general fact of discrete time processes that the Markov property implies the strong Markov property. That is for any stopping time $\tau$ with respect to $(\mc{F}_{t})_{t \geq 1}$, we have
\begin{equs}
    \pr(\sigma_{\tau+1}=\pi| \mathcal{F}_{\tau})=\pr(\sigma_{\tau+1}=\pi|\sigma_{\tau}),
\end{equs}
and 
\begin{equs}
    \pr(\sigma_{\tau+1}=\pi_1|\sigma_{\tau}=\pi_0) = \pr(\sigma_{1}=\pi_1| \sigma_{0}=\pi_0).
\end{equs}

Next we study Markovian properties of $\mrm{WP}(\lambda)$ with respect to a new filtration, $(\Lambda_t)_{t \geq 0}$. $(\Lambda_t)_{t \geq 0}$ will be defined by
\begin{equs}
    \Lambda_t:=\sigma(\lambda_0,\lambda_1,\dots,\lambda_t),\text{ where }\lambda_i = \mrm{shape}(\sigma_j).
\end{equs}
In words $\Lambda_t$ is the sigma algebra where we just keep track of the cycle information of the Weingarten process up to time $t$. 

It is important to note that $\mrm{WP}(\lambda)$ is in general, \textit{not a homogeneous Markov process with respect to $(\Lambda_k)_{k \geq 1}$}. To understand why this is the case, suppose at step $k+1$, $\sigma_{k+1} \in S_{n_{k+1}}$; then understanding the history of the cycle structures $(\lambda_j)_{j \leq k}$ give us more information about which cycle of $\sigma_{k+1}$ contains $n_{k+1}$ than just knowing $\lambda_k$. As a concrete example suppose we know the full history,
\begin{equs}
    (\lambda_0,\lambda_1)= ((3,1), (2,1,1)).
\end{equs}
Then there is a $1/3$ chance $\sigma_1(4)=4$ in which case $\sigma_2 \in S_3$, while if we just know $\lambda_1=(2,1,1)$, there is a chance $\lambda_0 \vdash 5$ in which case the conditional probability that $\sigma_1(4)=4$ is $1/2$. So the probability $\sigma_2 \in S_3$ conditioned only on $\lambda_1=(2,1,1)$ is larger than if we conditioned on the full history.  

Despite the preceding discussion, $\mrm{WP}(\lambda)$ has a Markovian property with respect to the filtration $(\Lambda_t)_{t \geq 1}$, if we restrict to \textit{pivotal times} which we define below

\begin{definition}\label{Pivot Change ST Def}
    For the process $\mathrm{WP}(\lambda)$, we define the pivot $n_t \in [n]$ of $\sigma_t$ to be the number such that $\sigma_t \in S_{n_t}$, and we define the $m$th \textit{pivotal time} to be, $\tau_m:= \min\{t: n_t = n-m\}$, that is $\tau_m$ tracks when the Weingarten process has changed level in the Weingarten graph. 
\end{definition}

\begin{remark}
    Note that $\tau_m$ is a stopping time with respect to the filtration $(\Lambda_t))_{t \geq 1}$. This is because we can figure out whether the level of the process has changed $m$ times by only keeping track of cycle structure up to a certain time.
\end{remark}
We now state a Markov-like property for $(\mrm{WP}(\lambda),(\Lambda_{t})_{t\geq 1}$, at the pivotal times whose proof we postpone until the end of the subsection.

\begin{lemma}\label{lm:New Level Markov Property}
    For $0\leq m \leq n$, and any permutation $\pi \in S_{n-m}$, then with respect to the process $\mrm{WP}(\lambda)$ for $\sigma \in S_n$ with $\mrm{shape}(\sigma)=\lambda$, we have,
    \begin{equs}
        \pr(\sigma_{\tau_m}=\pi | \Lambda_{\tau_m})=\pr(\sigma_{\tau_m}=\pi | \lambda_{\tau_m})=\frac{1}{|S_{\lambda_{\tau_m}}|}\mathbbm{1}_{\pi \in S_{\lambda_{\tau_m}}},
    \end{equs}
    and so as a corollary of Lemma \ref{Markov Property}, $(\sigma_{\tau_m},\sigma_{\tau_m+1},\dots,\sigma_{n+|\sigma|})$ conditioned on $\Lambda_{\tau_m}$ is distributed according to $\mrm{WP}(\lambda_{\tau_m})$. Additionally, for any stopping time $\tau$ which is almost surely a pivotal time (i.e. $\tau = \tau_m$ for some $m$),
    \begin{equs}\label{eq:general-pivotal-markov}
        \pr(\sigma_{\tau}=\pi | \Lambda_{\tau})=\pr(\sigma_{\tau}=\pi | \lambda_{\tau})=\frac{1}{|S_{\lambda_{\tau}}|}\mathbbm{1}_{\pi \in S_{\lambda_{\tau}}},
    \end{equs}
    and, $(\sigma_{\tau},\sigma_{\tau+1},\dots,\sigma_{n+|\sigma|})$ conditioned on $\Lambda_{\tau}$ is distributed according to $\mrm{WP}(\lambda_{\tau})$.
\end{lemma}

In order to prove the above Markov property, we will need the following lemma which gives us a bijection between sets of partial paths in the Weingarten graph with the same cycle types.

\begin{definition}
    Fix $0\leq m\leq  n$, $t \geq 0$, $\pi \in S_{n-m}$, and a sequence of partitions $(\eta_0,\dots,\eta_t)$ such that $\eta_t=\mrm{shape}(\pi)$. Let $\mrm{P}(\eta_0,\eta_1,\dots,\eta_{t};\pi)$ denote the set of directed paths $(\pi_0,\pi_1,\dots,\pi_t)$ in the unitary Weingarten graph taking $t$ edges, such that $\mrm{shape}(\pi_i)=\eta_i$ for $i=0,1,\dots,t$, ending at $\pi_t=\pi$, and such that the last edge from $\pi_{t-1}$ to $\pi_t$ is a dashed edge. 
\end{definition}

\begin{lemma}\label{lm:Partial Path Bijection}
    Fix a sequence of partitions $(\eta_0,\dots,\eta_t)$ with $\eta_t \vdash n-m \in [n-1]$, and two permutations $\pi,\pi' \in S_{\eta_t} \in S_{n-m}$, which are related by the conjugation $\pi'=\alpha\pi\alpha^{-1}$ for some $\alpha \in S_{n-m}$. Then the map $\mrm{P}(\eta_0,\eta_1,\dots,\eta_{t};\pi) \to\mrm{P}(\eta_0,\eta_1,\dots,\eta_{t};\pi')$ given by
    \begin{equs}
        (\pi_0,\dots,\pi_{t}) \mapsto (\alpha\pi_0\alpha^{-1},\dots,\alpha\pi_{t}\alpha^{-1}),
    \end{equs}
    is a bijection.
\end{lemma}
\begin{remark}
    As a clarification, when we write $\alpha\pi_j\alpha^{-1}$ where $\alpha \in S_{n-m}$ and $\pi_j \in S_{n'}$, for $n-m \leq n' \leq n$, we think of the product as an $S_{n'}$ permutation by identifying $\alpha$ with the permutation in $S_{n'}$ which fixes $\{n-m+1,\dots,n'-1,n'\}$.
    %whenever we write a product of permutations $\sigma_1 \sigma_2 \cdots \sigma_k$ where $\sigma_i \in S_{n_i}$ for possibly different values of $n_i$, we interpret this product belonging to $S_n$ for $n:=\max_i n_i$ as the product $\sigma_1'\sigma_2'\cdots \sigma_k'$ where each $\sigma_i' \in S_n$ is defined as acting identically to $\sigma_i$ on $\{1,2,\dots,n_i\}$ and as the identity on $\{n_i+1,\dots,n\}$.
\end{remark}

\begin{proof}
    Since every permutation in $S_{n-m}$ can be written as a product of transpositions, it suffices to prove the claim for $\alpha=(i\, j) \in S_{n-m}$ a transposition. The fact that this map and its inverse are one to one is trivial since for any $\rho,\rho' \in S_n$ with $\alpha \rho\alpha^{-1}=\alpha \rho' \alpha^{-1}$, it necessarily follows that $\rho = \rho'$. So the main part of this claim is actually that the map is well defined in the sense that it maps $\mrm{P}(\eta_0,\eta_1,\dots,\eta_{t}, ;\pi)$ to $\mrm{P}(\eta_0,\eta_1,\dots,\eta_{t}; \pi')$. 
    
    To see this, take $(\pi_0,\pi_1,\dots,\pi_t) \in \mrm{P}(\eta_0,\eta_1,\dots,\eta_{t}, \pi)$, and for any $ s\in \{0,1,\dots,t\} $ we must have $\pi_{s} \in S_{n'}$ for some $n-m \leq n' \leq n$. Then since $\alpha=(i \, j) \in S_{n-m} \subseteq S_{n'}$, we still have $\alpha \pi_{k'} \alpha^{-1} \in S_{n'}$, so $\alpha \pi_{s} \alpha^{-1}$ belongs to the same level of the Weingarten graph. Moreover, if $s<t$ and $\pi_{s}\in S_{n'}$, then as the last edge of the path $(\pi_0,\dots,\pi_t)$ is dashed, we must have the strict inequality $n'>n-m$. So if $s<t$ and $\pi_{s}$ is connected to $\pi_{s+1}$ by a directed dashed edge in the Weingarten graph, then $\pi_{s}$ fixes $n'$ and similarly so does $\alpha\pi_{s}\alpha^{-1}$ as $\alpha \in S_{n-m}$ and $n'>n-m$. As a result, $\alpha\pi_s\alpha^{-1}$ is connected to $\alpha\pi_{s+1}\alpha^{-1}$ by a directed dashed edge. Finally, if $\pi_{s}$ is connected to $\pi_{s+1}$ by a directed solid edge, then there exists a transposition $(r \, n') \in S_{n'}$ such that $\pi_{s+1}=\pi_s (r \, n')$, and thus
    \begin{equs}
        ~&\alpha \pi_{s+1} \alpha^{-1} = (i \, j) \pi_s (r \, n') (i \, j)\\
        & =\begin{cases}
            (i \, j) \pi_s  (i \, j)(r \, n') \hspace{5mm} \textit{if } r \notin \{i,j\},\\
            (i \, j)\pi_s (i \,j)(j \, n') \hspace{5mm} \textit{if } r= i, \\
            (i \, j)\pi_s (i \,j)(i \, n') \hspace{5mm} \textit{if } r= j,
        \end{cases}
    \end{equs}
    where we again use the fact that $n'>n-m$.
    In other words we still have $\alpha \pi_s \alpha^{-1}$ is connected to $\alpha \pi_{s+1}\alpha^{-1}$ by a solid directed edge in the Weingarten graph. This proves that $(\alpha \pi_0\alpha^{-1},\dots,\alpha \pi_{t}\alpha^{-1}) \in \mrm{P}(\eta_0,\eta_1,\dots,\eta_{t},\pi')$.
\end{proof}

\begin{proof}[Proof of Lemma \ref{lm:New Level Markov Property}]
    For any permutation $\pi \in S_{n-m}$ with $\pi \notin S_{\lambda_{\tau_m}}$, we must have $ \pr(\sigma_{\tau_m}=\pi | \Lambda_{\tau_m})=0$ as $\mrm{shape}(\sigma_{\tau_m})=\lambda_{\tau_m}$ by definition. Now for two distinct $\pi,\pi' \in S_{\lambda_{\tau_m}}$, there must exist a permutation $\alpha \in S_{n-m}$ with $\pi' = \alpha \pi \alpha^{-1}$. So for any particular $t \in \N$, and partitions $\eta_0,\eta_1,\dots,\eta_t$ of non-increasing size with $\eta_0 \vdash n$, $\eta_t \vdash n-m$, by Bayes' rule for conditional probabilities, we have
    \begin{equs}
        ~&\pr(\sigma_{\tau_m}=\pi|\{\tau_m = t\} \cap \{(\lambda_0,\dots,\lambda_t)=(\eta_0,\dots,\eta_t)\})\\
        &=\frac{\pr(\sigma_{\tau_m}=\pi \cap \{\tau_m = t\} \cap \{(\lambda_0,\dots,\lambda_t)=(\eta_0,\dots,\eta_t)\})}{\pr(\{\tau_m = t\} \cap \{(\lambda_0,\dots,\lambda_t)=(\eta_0,\dots,\eta_t)\})}\\
        &=\sum_{(\pi_0,\dots,\pi_{t-1},\pi_t) \in \mrm{P}(\eta_0,\dots,\eta_{t};\pi)} \frac{\pr((\sigma_0,\dots,\sigma_{t-1},\sigma_t)=(\pi_0,\dots,\pi_{t-1},\pi))}{\pr(\{\tau_m = t\} \cap \{(\lambda_0,\dots,\lambda_k)=(\eta_0,\dots,\eta_t)\})}\\
        &=\sum_{(\pi_0,\dots,\pi_{t-1},\pi_t) \in \mrm{P}(\eta_0,\dots,\eta_{t};\pi)}  \frac{\pr( (\alpha\sigma_0\alpha^{-1},\dots,\alpha\sigma_{t-1}\alpha^{-1},\alpha\sigma_t\alpha^{-1})=(\alpha\pi_0\alpha^{-1},\dots,\alpha\pi_{t-1}\alpha^{-1},\pi'))}{\pr( (\lambda_0,\dots,\lambda_t)=(\eta_0,\dots,\eta_t))}\\   &=\sum_{(\pi_0,\dots,\pi_{t-1},\pi_t) \in \mrm{P}(\eta_0,\dots,\eta_{t};\pi')}  \frac{\pr((\sigma_0,\dots,\sigma_{t-1},\sigma_t)=(\pi_0,\dots,\pi_{t-1},\pi'))}{\pr(\{\tau_m = t\} \cap \{(\lambda_0,\dots,\lambda_t)=(\eta_0,\dots,\eta_t)\})}\\
        &=\pr(\sigma_{\tau_m}=\pi'|\{\tau_m = t\} \cap \{(\lambda_0,\dots,\lambda_t)=(\eta_0,\dots,\eta_t)\}),
    \end{equs}
    where in the second to last equality, we applied the bijection of Lemma \ref{lm:Partial Path Bijection}, and the fact that the probability of $(\sigma_0,\dots,\sigma_t)=(\pi_0,\dots,\pi_t)$ for any fixed $(\pi_0,\dots,\pi_t)$ is only a function of $(\mrm{shape}(\pi_0),\dots,\mrm{shape}(\pi_t))$; in particular, this probability can be expressed as a ratio of products of Catalan numbers of cycle lengths. So we have proven that $\sigma_{\tau_m}$ is uniform over $S_{\lambda_{\tau_m}}$, conditional on $\Lambda_{\tau_m}$. Or in other words,
    \begin{equs}
        \pr(\sigma_{\tau_m}=\pi | \Lambda_{\tau_m})=\frac{1}{|S_{\lambda_{\tau_m}}|}\mathbbm{1}_{\pi \in S_{\lambda_{\tau_m}}}.
    \end{equs}
     It follows at once from the tower property of conditional expectations that we also have,
     \begin{equs}
         \pr(\sigma_{\tau_m}=\pi | \lambda_{\tau_m})=\frac{1}{|S_{\lambda_{\tau_m}}|}\mathbbm{1}_{\pi \in S_{\lambda_{\tau_m}}}.
     \end{equs}
     The proof of \eqref{eq:general-pivotal-markov} is identical so we leave it out.
\end{proof}

\subsection{Outline of the  Proof of Theorem \ref{Main Weingarten Bound}.}

In order to prove Theorem \ref{Main Weingarten Bound}, it suffices to prove that $\frac{|\mathrm{P}(\sigma,|\sigma|+2g)|}{|\mathrm{P}(\sigma,|\sigma|)|} \leq (Cn^3)^g$ for $C = 6\sqrt{8}\times 10^6$. Moreover, by a simple inductive argument we will only need to show $\frac{|\mathrm{P}(\sigma,|\sigma|+2)|}{|\mathrm{P}(\sigma,|\sigma|)|} \leq Cn^3$ for some $C>0$. To prove this, we start by closely following the arguments of \cite{Collins2017} by observing that any path $ \mathrm{p} \in \mathrm{P}(\sigma,|\sigma|+2)$ has exactly one ``defect'' solid edge $\sigma_j \to \sigma_{j+1}$ where $|\sigma_{j+1}|>|\sigma_j| $. If we count all the paths in $ \mathrm{p} \in \mathrm{P}(\sigma,|\sigma|+2)$ by counting all possible ``partial paths'' until the defect, and after the defect, followed by a worst case estimate involving the ratio of Catalan numbers capturing the effect of the defect, we are led to the bound,
\begin{equs}
    \frac{|\mathrm{P}(\sigma,|\sigma|+2)|}{|\mathrm{P}(\sigma,|\sigma|)|} &\leq 6\sqrt{8}n \sum_{j=0}^{n+|\sigma|}\E_{\mathrm{Uniform}(\mathrm{P}(\sigma,|\sigma|))}[L_j^{3/2}]\\
    &=6\sqrt{8}n \sum_{j=0}^{n+|\sigma|}\E_{\mathrm{WP}(\sigma)}[L_j^{3/2}],
\end{equs}
where $L_j$ denotes the length of the longest cycle of $\sigma_j$ (the $j$th permutation in the path). The key insight is then observing that the path counts are class functions on $S_n$, so we can take the average over all $\sigma' \in S_{\lambda}$ for $\lambda=\mrm{shape}(\sigma)$ to obtain,
\begin{equs}
    \frac{|\mathrm{P}(\sigma,|\sigma|+2)|}{|\mathrm{P}(\sigma,|\sigma|)|} \leq 6\sqrt{8}n \sum_{j=0}^{n+|\sigma|}\E_{\mathrm{WP}(\lambda)}[L_j^{3/2}],
\end{equs}
 The proof of this formula is the content of Subsection \ref{Reduction to Expectation Section}.

Next in Subsection \ref{Bounds on expectation section}, we finish the argument by proving the bound, 
\begin{equs}
    \sum_{j=0}^{n+|\sigma|}\E_{\mathrm{WP}(\lambda)}[L_j^{3/2}] \leq 10^6 n^{3/2}L_0^{1/2} \leq 10^6 n^2.
\end{equs}
The key insight will be, Lemma \ref{lm:time-to-halve}, that every $O(L_0^{1/2}(1+\log(n/L_0)^2))$ steps, $L_t$ shrinks by a factor $2/3$. This step critically relies on the extra randomness of working with $\mrm{WP}(\lambda)$ instead of $\mrm{WP}(\sigma)$. The idea here is to treat every cycle at least $2/3$ the size of the largest one at the start as a coupon which gets collected as soon as the cycle shrinks below the threshold of $2/3$ the size of the largest cycle at the start. We will find that if there are $k$ coupons, the wait time between collecting the $i$th and $(i+1)$th coupon will be approximately a geometric random variable with parameter $p \approx\frac{L_0^{1/2}}{n}(k-i)$. This comparison to a geometric random variable is the content of Lemma \ref{lm:T^{(i)}}, and is the most involved step of the proof. It will be broken into four steps, which at a very high level are summarized as follows: Step 1 is to show that pivotal times happen frequently, step 2 is to lower bound the probability that we encounter a large cycle at a pivotal time, step 3 is to lower bound the probability that a large pivotal cycle will split into two similar
sized pieces in a given step, and step 4 is to show that a pivotal time will likely occur shortly after this large cycle gets broken up. 

Once Lemma \ref{lm:T^{(i)}} is established, the remainder of the proof of Lemma \ref{lm:time-to-halve} proceeds similarly to estimates for the usual coupon collector problem. The proof here is not completely optimal (there's an extra unnecessary logarithm), but suffices for our purpose.

\subsection{Estimating $\mathrm{P}(\sigma, |\sigma|+2g)$.}\label{Reduction to Expectation Section}

We start by recalling several definitions from \cite{Collins2017}. For a fixed path $p \in \mathrm{P}(\sigma,n+l)$ we can define $j:= j(p)$ such that the $(j+1)$th solid arrow crossed in the path is the first transition $\sigma_{j+n-r} \to \sigma_{j+n-r+1}$ such that $|\sigma_{j+n-r}|+1 = |\sigma_{j+n-r+1}|$, and we define
\begin{equs}
    \mathrm{P}_j(\sigma,l):=\{p \in \mathrm{P}(\sigma,l): j(p)=j\}.
\end{equs}
We clearly have
\begin{equs}\label{First Backtack Decomp}
    \mathrm{P}(\sigma,l)= \bigcup_{j=0}^{|\sigma|} \mathrm{P}_j(\sigma,l). 
\end{equs}

Moreover, every path can be decomposed into the path up to this step, and the remaining part, i.e. for $|\sigma|=k$,
\begin{equs}\label{First Backtack, level Decomp}
    \mathrm{P}_j(\sigma,l) = \bigcup_{r=2}^{k} \bigcup_{\rho \in S_r} \bigcup_{\substack{\rho' \in S_r \\ |\rho'|=|\rho|+1\\
    \rho \to \rho'}} \mathrm{P}_j(\sigma,\rho,\rho',l),
\end{equs}
where $\mathrm{P}_j(\sigma,\rho,\rho',l)$ denotes the set of paths $\mathrm{P}_j(\sigma,l)$, where the first edge $\sigma_{j+n-r} \to \sigma_{j+n-r+1}$ with $|\sigma_{j+n-r}|+1 = |\sigma_{j+n-r+1}|$ occurs at $\sigma_{j+n-r} =\rho$ and $\sigma_{j+k-r+1}=\rho'$.

\noindent If we let $\Tilde{\mathrm{P}}_j(\sigma,\rho)$ denotes the set of all "partial paths" from $\sigma \to \rho$, i.e. paths $\sigma= \sigma_0, \sigma_1,...,\sigma_{j+n-r}=\rho$ with no increasing norm steps. Then we have the factorization,

\begin{equs}\label{Partial path decomposition}
    \mathrm{P}_j(\sigma,\rho,\rho',|\sigma|+2g) \simeq \Tilde{\mathrm{P}}_j(\sigma,\rho) \times \mathrm{P}(\rho', |\rho'|+2g-2).
\end{equs}
Next recall the following bound from \cite{Collins2017}.

\begin{lemma}\cite[Lemma 3.4]{Collins2017}\label{Catalan quotient bound}
    For any $j,k \in \N$ with $j < k$,
    \begin{equs}
        \max_{1 \leq j \leq k} \frac{\mathrm{Cat}(k-1)}{\mathrm{Cat}(j-1)\mathrm{Cat}(k-j-1)} \leq 6 k^{3/2}.
    \end{equs}
\end{lemma}

\begin{remark}
    This Catalan quotient bound is not the statement of \cite[Lemma 3.4]{Collins2017}, but it is in the proof over the course of the proof of the lemma.
\end{remark}

Finally we need one last definition before proving the main lemma of this section.
\begin{definition}
    For the Weingarten process $\mathrm{WP}(\lambda)$, we let $L_j$ be the random variable taking as its value, the length of the longest cycle in the permutation $\sigma_j$ (the $j$th permutation in the Weingarten process), and we set $L_j=0$ if $j>n+|\sigma|$.
\end{definition}

The main lemma of this section is the following bound comparing the first order coefficient to the zeroth order coefficient in the $1/N$ expansion of the Weingarten function. The proof of this lemma follows closely the proof of \cite[Theorem 3.1]{Collins2017}, in particularly sections 3.3.3 and 3.3.4. The main difference is only at the end, when we apply Lemma \ref{Unitary Process is uniform over paths}, nevertheless we give the details for the sake of completeness.

\begin{lemma}\label{Single Defect Path Bound}
    For any $\sigma \in S_n$ with cycle type $\mrm{shape}(\sigma)=\lambda$
    \begin{equs}\label{Path enumeration bounded by expecations}
        \frac{|\mathrm{P}(\sigma,|\sigma|+2)|}{|\mathrm{P}(\sigma,|\sigma|)|} \leq 6\sqrt{8}n\sum_{i=0}^{n+|\sigma|}\E_{\mathrm{WP}(\lambda)}[L_i^{3/2}].
    \end{equs}
\end{lemma}

\begin{remark}
    The proof of this estimate is almost exactly the same as the proof of \cite[Theorem 3.1]{Collins2017}, except that they bound the length of the longest cycle of a permutation in $S_{n-r}$ by $S_n$, while we keep the right hand-side of our bound in terms of the longest cycle length.
\end{remark}

\begin{proof}

First using the path decomposition \eqref{First Backtack Decomp}, \eqref{First Backtack, level Decomp}, and \eqref{Partial path decomposition} successively we have
\begin{equs}
    |\mathrm{P}(\sigma,|\sigma|+2)| &= \sum_{j=0}^{|\sigma|} |\mathrm{P}_j(\sigma, |\sigma|+2)|\\
    &= \sum_{j=0}^{|\sigma|} \sum_{r=2}^{n} \sum_{\rho \in S_r} \sum_{\substack{\rho' \in S_r \\ |\rho'|=|\rho|+1\\
    \rho \to \rho'}} |\mathrm{P}_j(\sigma,\rho,\rho' |\sigma|+2)| \\
    &= \sum_{j=0}^{|\sigma|} \sum_{r=2}^{n} \sum_{\rho \in S_r} \sum_{\substack{\rho' \in S_r \\ |\rho'|=|\rho|+1\\
    \rho \to \rho'}}|\Tilde{\mathrm{P}}_j(\sigma,\rho)| |\mathrm{P}(\rho',|\rho'|)|.  
\end{equs}

We can now rewrite $|\mathrm{P}(\rho',|\rho'|)|=|\mathrm{P}(\rho,|\rho|)|\cdot \frac{|\mathrm{P}(\rho',|\rho'|)|}{|\mathrm{P}(\rho,|\rho|)|}$ and bound $\frac{|\mathrm{P}(\rho',|\rho'|)|}{|\mathrm{P}(\rho,|\rho|)|}$ by its maximum value over all choices of $\rho'$,
\begin{equs}
    |\mathrm{P}(\sigma,|\sigma|+2)|&=\sum_{j=0}^{|\sigma|} \sum_{r=2}^{n} \sum_{\rho \in S_r} \sum_{\substack{\rho' \in S_r \\ |\rho'|=|\rho|+1\\
    \rho \to \rho'}}|\Tilde{\mathrm{P}}_j(\sigma,\rho)| |\mathrm{P}(\rho,|\rho|)|\cdot \frac{|\mathrm{P}(\rho',|\rho'|)|}{|\mathrm{P}(\rho,|\rho|)|}\\
    &\leq \sum_{j=0}^{|\sigma|} \sum_{r=2}^{n} \sum_{\rho \in S_r} \sum_{\rho' \in S_r, |\rho'|=|\rho|+1,\rho \to \rho'} |\Tilde{\mathrm{P}}_j(\sigma,\rho)| |\mathrm{P}(\rho,|\rho|)| \max_{\rho'' \to \rho} \frac{|\mathrm{P}(\rho'',|\rho''|)|}{|\mathrm{P}(\rho,|\rho|)|}.
\end{equs}
Now, let $LC(\rho)$ denote the size of the longest cycle in a permutation $\rho$. Then plugging in the enumeration formula Lemma \ref{Path Enumeration Formula}, $\frac{|\mathrm{P}(\rho'',|\rho''|)|}{|\mathrm{P}(\rho,|\rho|)|}. = \frac{\mrm{Cat}(k-1)}{\mrm{Cat}(j-1)\mrm{Cat}(k-j-1)}$ for some $1 \leq j \leq k \leq 2\mrm{LC}(\rho)$. So applying the Catalan number quotient bound Lemma \ref{Catalan quotient bound}, we have
\begin{equs}
     |\mathrm{P}(\sigma,|\sigma|+2)|&\leq \sum_{j=0}^{|\sigma|} \sum_{r=2}^{n} \sum_{\rho \in S_r} \sum_{\rho' \in S_r, |\rho'|=|\rho|+1,\rho \to \rho'} |\Tilde{\mathrm{P}}_j(\sigma,\rho)| |\mathrm{P}(\rho,|\rho|)| \cdot 6(2LC(\rho))^{3/2}\\
    &\leq 6\sqrt{8}n \sum_{j=0}^{|\sigma|} \sum_{r=2}^{n} \sum_{\rho \in S_r}  |\tilde{\mathrm{P}}_j(\sigma,\rho)| |\mathrm{P}(\rho,|\rho|)| LC(\rho)^{3/2},
\end{equs}
where in the last inequality, we used the fact that there were at most $r-1 \leq n$ transpositions $\rho'$ to sum over.

At this point observe that there is a natural mapping 
\begin{equs}
    \bigcup_{j=0}^{|\sigma|}\bigcup_{r=2}^{n}\bigcup_{\rho \in S_r}\tilde{\mathrm{P}}_j(\sigma,\rho)\times \mathrm{P}(\rho,|\rho|) \to \mrm{P}(\sigma,|\sigma|)\times \{0,1,\dots,n+|\sigma|\},
\end{equs}
mapping
\begin{equs}
    (\sigma_0,\dots,\sigma_{n-r+j-1},\rho) \times (\rho,\sigma_1',\dots, \sigma_{|\rho|+r-1}, \varnothing) \mapsto (\sigma_0,\dots,\sigma_{n-r+j-1},\rho,\sigma_1',\dots, \sigma_{|\rho|+r-1}, \varnothing) \times \{n-r+j\}.
\end{equs}
Moreover this map is obviously injective as the number in $\{0,1,\dots,n+|\sigma|\}$ uniquely determines how to invert the map by tracking where to split the path. Thus,

\begin{equs}
    |\mathrm{P}(\sigma,|\sigma|+2)|&\leq 6\sqrt{8}n \sum_{j=0}^{|\sigma|} \sum_{r=2}^{n} \sum_{\rho \in S_r}  |\Tilde{\mathrm{P}}_j(\sigma,\rho)| |\mathrm{P}(\rho,|\rho|)| LC(\rho)^{3/2}\\
    &\leq 6\sqrt{8}n 
    \sum_{(\sigma_0,\sigma_1,\dots,\sigma_{n+|\sigma|}) \in \mrm{P}(\sigma,|\sigma|)} \sum_{i=0}^{n+|\sigma|}LC(\sigma_i)^{3/2}\\
    & = 6\sqrt{8}n |\mathrm{P}(\sigma,|\sigma|)|\sum_{i=0}^{n+|\sigma|}\E_{\mathrm{Unif}(\mathrm{P}(\sigma,|\sigma|))}[LC(\sigma_i)^{3/2}]\\
    &=6\sqrt{8}n |\mathrm{P}(\sigma,|\sigma|)|\sum_{i=0}^{n+|\sigma|}\E_{\mathrm{WP}(\sigma)}[LC(\sigma_i)^{3/2}],
\end{equs}
where we applied Lemma \ref{Unitary Process is uniform over paths} in the last equality. Finally we now recall that the Weingarten function is a class function, as it can be expressed purely in terms of characters of the symmetric group, since $\mrm{P}(\sigma,|\sigma|+2)$ must also be a class function as it appears in the $\frac{1}{N}$ expansion of $\mrm{Wg}^{\mrm{U}(N)}_N(\sigma)$. Thus, letting $\lambda = \mrm{shape}(\sigma)$,
\begin{equs}
    |\mathrm{P}(\sigma,|\sigma|+2)| &= \frac{1}{|S_{\lambda}|}\sum_{\sigma' \in S_{\lambda}}|\mathrm{P}(\sigma',|\sigma'|+2)|\\
    &\leq \frac{1}{|S_{\lambda}|}\sum_{\sigma' \in S_{\lambda}} 6\sqrt{8}n |\mathrm{P}(\sigma',|\sigma'|)|\sum_{i=0}^{n+|\sigma'|}\E_{\mathrm{WP}(\sigma')}[LC(\sigma_i)^{3/2}]\\
    &=6\sqrt{8}n |\mathrm{P}(\sigma,|\sigma|)|\sum_{i=0}^{n+|\sigma|}\E_{\mathrm{WP}(\lambda)}[LC(\sigma_i)^{3/2}].
\end{equs}
\end{proof}
From now on, we assume all probabilities and expectations are taken with respect to $\mathrm{WP}(\lambda)$ unless explicitly stated otherwise. It now suffices to prove the following proposition.

\begin{prop}\label{prop:Main expectation bound}
    For any $n \geq 1$ and any $\sigma \in S_n$ with longest cycle $L_0$,
    \begin{equs}
        \sum_{i=0}^{n+|\sigma|}\E[L_j^{3/2}] \leq 10^6 L_0^{1/2}n^{3/2} \leq 10^6 n^2.
    \end{equs}
\end{prop}

The proof of this bound will be deferred to Section \ref{Bounds on expectation section}. In the remainder of the section we assume Proposition \ref{prop:Main expectation bound}, and complete the proof of Theorem \ref{Main Weingarten Bound}.

\begin{cor}\label{Single Defect Lemma}
    For all $n$ and $\sigma \in S_n$,
    \begin{equs}
        \frac{|\mathrm{P}(\sigma,|\sigma|+2)|}{|\mathrm{P}(\sigma,|\sigma|)|} \leq 6\sqrt{8}\times 10^6 \cdot n^{3}.
    \end{equs}
\end{cor}

\begin{proof}
    Plug in the inequality of Proposition \ref{prop:Main expectation bound} into Lemma \ref{Single Defect Path Bound}.
\end{proof}

\begin{lemma}
    For all $n$ and $\sigma \in S_n$,
    \begin{equs}\label{Unitary Genus g bound}
        \frac{|\mathrm{P}(\sigma,|\sigma|+2g)|}{|\mathrm{P}(\sigma,|\sigma|)|} \leq (6\sqrt{8}\times 10^6 \times n^{3} )^g.
    \end{equs}
\end{lemma}

\begin{proof}
    The case $g=0$ requires no proof, and the case $g=1$ is the statement of the previous lemma. The general case now follows by induction in a similar manner to the proof of Lemma \ref{Single Defect Path Bound}, for notational convenience setting $C=6\sqrt{8}\times 10^6$,
    \begin{equs}
        |\mathrm{P}(\sigma,|\sigma|+2g)| &= \sum_{j=0}^{|\sigma|} |\mathrm{P}_j(\sigma, |\sigma|+2g)|\\
        &= \sum_{j=0}^{|\sigma|} \sum_{r=2}^{k} \sum_{\rho \in S_r} \sum_{\rho' \in S_r, |\rho'|=|\rho|+1,\rho \to \rho'} |\mathrm{P}_j(\sigma,\rho,\rho' |\sigma|+2g)| \\
        &= \sum_{j=0}^{|\sigma|} \sum_{r=2}^{k} \sum_{\rho \in S_r} \sum_{\rho' \in S_r, |\rho'|=|\rho|+1,\rho \to \rho'}|\Tilde{\mathrm{P}}_j(\sigma,\rho)| |\mathrm{P}(\rho',|\rho'|+2g-2)|  \\
        &\leq (Cn^3)^{g-1}\sum_{j=0}^{|\sigma|} \sum_{r=2}^{k} \sum_{\rho \in S_r} \sum_{\rho' \in S_r, |\rho'|=|\rho|+1,\rho \to \rho'} |\Tilde{\mathrm{P}}_j(\sigma,\rho)| |\mathrm{P}(\rho',|\rho'|)|\\
        & \leq (Cn^3)^{g-1} |\mathrm{P}(\sigma,|\sigma|+2)|\\
        &\leq (Cn^3)^{g} |\mathrm{P}(\sigma,|\sigma|)|.
    \end{equs}
\end{proof}

\noindent
\begin{proof}{Proof of Theorem \ref{Main Weingarten Bound}:} 
Directly plugging in the bound \eqref{Unitary Genus g bound} into the path expansion \eqref{Path Expansion Unitary}, we have
\begin{equs}
    |N^{-n-|\sigma|}\mathrm{Wg}_N^{\mathrm{U}}(\sigma)|&\leq |\mathrm{P}(\sigma,|\sigma|)|\sum_{g \geq 0} (Cn^3)^g N^{-2g}\\
    &= \frac{|\mathrm{P}(\sigma,|\sigma|)|}{1-Cn^3/N^2},
\end{equs}
for $Cn^3<N^2$. Now dividing both sides by $|\mathrm{P}(\sigma,|\sigma|)|$ and recalling that $|\mathrm{P}(\sigma,|\sigma|)|=|\mathrm{Moeb}(\sigma)|$, we get the desired bound.

\end{proof}

\subsection{Proof of Proposition \ref{prop:Main expectation bound}}\label{Bounds on expectation section}

Recall from Definition \ref{Pivot Change ST Def}, that $\tau_m$ tracks the first time in $\mathrm{WP}(\lambda)$, when the permutation in the Weingarten process is an element of $S_{n-m}$, and that $\tau_m$ is a stopping time with respect to $(\mc{F}_t)_{t \geq 1}$ and $(\Lambda_t)_{t\geq 1}$. We now introduce a few more random processes.

\begin{definition}
    Define $T = \tau_{\tilde{i}}$ where $\tilde{i} = \min\{i > 0: L_{\tau_{i}}<\frac{2}{3}L_{0}\}$.
\end{definition}
It turns out Proposition \ref{prop:Main expectation bound} follows rather directly by having the correct bound on the expected time until the largest cycle of $\sigma$ shrinks to $2/3$ the size, $\E[T]$.

\begin{lemma}\label{lm:time-to-halve}
    Suppose $n> 6$. For any $\sigma \in S_n$ with longest cycle $L_0 \geq 6$, 
    \begin{equs}
        \E[T] \leq 10^4\frac{n}{L_0^{1/2}}\bigg(1+\log\bigg(\frac{3n}{2L_0}\bigg)\bigg)^2.
    \end{equs}
    with respect to the Weingarten process $\mrm{WP}(\lambda)$.
\end{lemma}

\begin{proof}[Proof of Proposition \ref{prop:Main expectation bound} assuming Lemma \ref{lm:time-to-halve}] We apply strong induction on $n$. 

For the base case when $n \leq 6$ or when $L_0 < 6$, the claim is obvious and there is nothing to prove.

For the inductive step, we split the sum $\sum_{j=0}^{n+|\sigma|} \E[L_j^{3/2}]$ into two parts at the index $T$, bounding the $L_j$ until $T$ by $L_0$, and applying the bound of Lemma \ref{lm:time-to-halve} to the first part, while applying the Markov property, \eqref{eq:general-pivotal-markov}, and the inductive hypothesis to the second part as follows,
    \begin{equs}
        ~&\sum_{j=0}^{n+|\sigma|} \E[L_j^{3/2}] \\
        &=\E\bigg[\sum_{j=0}^{T-1} L_j^{3/2}\bigg] +\E\bigg[\sum_{j=T}^{n+|\sigma|} L_j^{3/2}\bigg]\\
        &\leq \E[T]L_0^{3/2}+\E_{\mrm{WP}(\lambda)}\bigg[\E_{\mrm{WP}(\lambda_T)}\bigg[\sum_{j=0}^{n_T+|\sigma_T|} L_j^{3/2}\bigg]\bigg]\\
        &\leq 10^4nL_0\bigg(1+\log\bigg(\frac{3n}{2L_0}\bigg)\bigg)^2+ \E_{\mrm{WP}(\lambda)}[10^6 n_T^{3/2}L_T^{1/2}]\\
        &\leq 10^4nL_0 \times 4 (3n/2L_0)^{1/2}+(2/3)^{1/2}10^6 n^{3/2}L_0^{1/2}\\
        &\leq 10^6 n^{3/2}L_0^{1/2}
    \end{equs}
    where we used the elementary inequalities $1+\log(x) \leq 4x^{1/4}$ for $x \geq 1$, and $4\times (3/2)^{1/2}+100(2/3)^{1/2}\leq 100$.
    %where we used the fact that $L_T \leq 2L_0/3$ by definition. Lastly the proposition follows from the inequality $n x\log(n/x) \leq n^2$ for $0<x \leq n$ which can be verified by elementary calculus.
\end{proof}

The rest of the subsection is dedicated to proving Lemma \ref{lm:time-to-halve}. Before diving into the proof, we need to define a few more stopping times. For the convenience of the reader, in Table \ref{table:RVs} we give a description of many of the random variables appearing throughout the rest of the section.
\begin{table}
    \centering
    \begin{tabular}{|c|c|}
        \hline
       Random Variable  &  Description\\
       \hline
         $n_t$ & Pivotal element \\
         \hline
         $C_t$ & Length of cycle containing pivotal element \\
         \hline
         $L_t$ & Length of longest cycle \\
         \hline
         $\tau_i$ & Times when pivotal element changes (pivotal times). \\
         \hline
         $t_i$ & Times when a cycle of size $\geq \frac{2}{3}L_0$ splits into two cycles of size $< \frac{2}{3}L_0$\\
         \hline
         $T_i$ & Next pivotal time after $t_i$\\
         \hline
         $T$ & First pivotal time when no cycles of size $\geq \frac{2}{3}L_0$ remain\\
         \hline
    \end{tabular}
    \caption{Random variables appearing in the proof of Lemma \ref{lm:time-to-halve}.}
    \label{table:RVs}
\end{table}

Suppose $\sigma_0$ has $k$ cycles of size $\geq \frac{2}{3}L_0$, for $i<k$,
\begin{equs}
    ~&t_0 = 0\\
    & t_{i} = \inf\{t \geq 0: \sigma_t \text{ has } k-i \text{ cycles of length } \geq \frac{2}{3}L_0\}.
\end{equs}
It is clear from the definition that the times $t_i$ are stopping times with respect to $(\Lambda_t)$ and $(\mc{F}_t)$, and that $t_0<t_1<t_2<\dots <t_k$. For the convenience of being able to apply the Markov property with respect to $(\Lambda_t)$, Lemma \ref{lm:New Level Markov Property}, we instead work with the following stopping times,

\begin{equs}
    T_{i} = \inf\{t \geq t_i: T_i \text{ is pivotal}\}.
\end{equs}
Recall that a time $t$ is pivotal if $t=\tau_m$ for some $m \leq n$. It is also true that the times $T_i$ are stopping times with respect to $(\Lambda_t)$ and $(\mc{F}_t)$, and that $T_0<T_1<T_2<\dots <T_k$. It is also obvious by definition that
\begin{equs}
    T=T_k,
\end{equs}
so bounding $\E[T]$ reduces to understanding the distribution of $T_{i+1}-T_i$ for $i<k$. The next core lemma explains roughly that the increments $T_{i+1}-T_i$ are distributed like geometric random variables with parameter $p \sim (k-i)\frac{L_0^{1/2}}{n}$. This means that with high probability, $T_{i+1}-T_i$ should not be much larger than $\frac{n}{L_0^{1/2}(k-i)}$ up to constants. The proof is divided into four steps which are summarized below.

\textbf{Step 1:} Every 2 steps the Weingarten process sees a pivotal cycle with $\geq 1/4$ probability. This is because every time a cycle of size $\ell>1$ splits, it splits into pieces of size $1$ and $\ell-1$ with $\geq 1/4$ probability due to a simple bound on the ratio of consecutive Catalan numbers.

\textbf{Step 2:} At a pivotal time between $T_i$ and $T_{i+1}$, the Weingarten process sees a large cycle of size $\geq \frac{2}{3}L_0$ as the pivotal cycle with probability at least $\frac{2(k-i)L_0}{3n}$. This step critically relies on the Markov property of $\mrm{WP}(\lambda)$ with respect to $\Lambda_t$ at pivotal times. The idea is that given its cycle type, the permutation at this pivotal time is uniform, so for any of the $k-i$ cycles of length $\geq \frac{2}{3}L_0$, there is a $\geq \frac{2L_0}{3n}$ chance it contains the pivotal element.

\textbf{Step 3:} In this step, we observe that if we are at a pivotal time, and the pivotal cycle has length $\geq \frac{2}{3}L_0$, then there is at least an $O(L_0^{-1/2})$ chance it splits into two cycle of length $< \frac{2}{3}L_0$. This step comes from directly bounding the probability that the pivotal cycle splits into two cycles of similar sizes. To bound this probability, we express it as a sum of quotients of Catalan numbers and apply Lemma \ref{Catalan quotient bound} to bound each of the Catalan number quotients individually.

\textbf{Step 4:} Combining steps 1-3, we see that $t_{i+1}$ has occurred within 3 steps after $t$ with probability $O((k-i)L_0^{1/2}n^{-1})$. Finally recycling the argument of step 1, it follows that with $\geq \frac{1}{4}$ probability, $T_{i+1}$ happens within two steps of $t_{i+1}$, concluding the argument.

Let us now give the precise statement and proof.

\begin{lemma}\label{lm:T^{(i)}}
    Fix any $t \in [n+|\sigma|-6]$, and suppose $L_0 \geq 6$. Then
    \begin{equs}
        \pr(T_{i+1}> t+5|T_{i}\leq t <T_{i+1}) \leq 1-10^{-3}\frac{(k-i)L_0^{1/2}}{n},
    \end{equs}
    and as a corollary,
    \begin{equs}\label{eq:T^i tail}
        \pr(T_{i+1}-T_i > 5t) \leq \bigg(1-10^{-3}\frac{(k-i)L_0^{1/2}}{n}\bigg)^t,
    \end{equs}
    with respect to $\mrm{WP}(\lambda)$
\end{lemma}

\begin{proof}[Proof of Lemma \ref{lm:T^{(i)}}]
    \textbf{Step 1:} We first prove that
    \begin{equs}\label{eq:Pivotal in two steps}
        \pr(t+1 \text{ or } t+2 \text{ is pivotal}|T_i\leq t < T_{i+1} ) \geq \frac{1}{4}.
    \end{equs}
    Recall a time $t$ is pivotal if it equals $\tau_m$ for some $m$. Since $\{T_i\leq t < T_{i+1}\}$ is an $\mc{F}_{t}$ measurable event, it will suffice to show that in general,
    \begin{equs}
        \pr(t+1 \text{ or } t+2 \text{ is pivotal}|\mc{F}_t) \geq \frac{1}{4}.
    \end{equs}
    Moreover by the Markov property, this reduces to showing that for any permutation $\pi$, 
    \begin{equs}
        \pr(t+1 \text{ or } t+2 \text{ is pivotal}|\sigma_t=\pi) \geq \frac{1}{4}.
    \end{equs} 
    Letting $\ell$ denote the length of the pivotal cycle in $\pi$. We split into cases. If $\ell=1$, then by definition $t+1$ will necessarily be pivotal so $\pr(t+t\text{ is pivotal}| \sigma_t=\pi)=1$. On the other hand $t+2$ will be pivotal iff the pivotal cycle at $t+1$ has length $1$. Thus, if $\ell>1$,
    \begin{equs}
        \pr(t+2\text{ is pivotal}| \sigma_t=\pi) = \frac{\mrm{Cat}(\ell-2)\mrm{Cat}(0)}{\mrm{Cat(\ell-1)}} \geq \frac{1}{4},
    \end{equs}
    where we have used the elementary inequality that the ratio of consecutive Catalan numbers is bounded above by $4$. 

    \textbf{Step 2:} At this point, define $C_t$ to be the length of the cycle containing the pivotal element $n_t$ at time $t$. Recall $n_t \in [n]$ is defined as the number such that $\sigma_t \in S_{n_t}$. The next step is to show that 
    \begin{equs}\label{eq:pivotal-cycle-appears}
        \pr(C_{t+2}\geq \frac{2}{3}L_0|T_i\leq t < T_{i+1})+\pr(C_{t+1}\geq \frac{2}{3}L_0|T_i\leq t < T_{i+1})\geq \frac{(k-i)L_0}{6n}.
    \end{equs}
    It will suffice to show that both 
    \begin{equs}\label{eq:See-big-cycle-at-pivot}
        ~&\pr(C_{t+1}\geq \frac{2}{3}L_0|t+1 \text{ is pivotal, } T_i\leq t < T_{i+1}) \geq \frac{2(k-i)L_0}{3n}, \text{ and}\\
        &\pr(C_{t+2}\geq \frac{2}{3}L_0|t+2 \text{ is pivotal, } T_i\leq t < T_{i+1}) \geq \frac{2(k-i)L_0}{3n},
    \end{equs}
    since if this were the case then,
    \begin{equs}
        ~&\pr(C_{t+2}\geq \frac{2}{3}L_0|T_i\leq t < T_{i+1})+\pr(C_{t+1}\geq \frac{2}{3}L_0|T_i\leq t < T_{i+1})\\
        &\geq \pr(C_{t+2} \geq \frac{2}{3}L_0, \, t+2 \text{ is pivotal}| T_i\leq t < T_{i+1})+\pr(C_{t+1} \geq \frac{2}{3}L_0,\, t+1 \text{ is pivotal}| T_i\leq t < T_{i+1})\\
        &=\pr(C_{t+2} \geq \frac{2}{3}L_0| t+2 \text{ is pivotal}, T_i\leq t < T_{i+1}) \pr(t+2 \text{ is pivotal}| T_i\leq t < T_{i+1})\\
        &+\pr(C_{t+1} \geq \frac{2}{3}L_0| t+1 \text{ is pivotal}, T_i\leq t < T_{i+1}) \pr(t+1 \text{ is pivotal}| T_i\leq t < T_{i+1})\\
        &\geq \frac{2(k-i)L_0}{3n}(\pr(t+2 \text{ is pivotal}| T_i\leq t < T_{i+1})+\pr(t+1 \text{ is pivotal}| T_i\leq t < T_{i+1}))\\
        &\geq \frac{(k-i)L_0}{6n},
    \end{equs}
    where in the last line, we applied \eqref{eq:Pivotal in two steps}.
    Now we only prove the first inequality in \eqref{eq:See-big-cycle-at-pivot} as the second is proven in an identical manner. We first apply Baye's formula for conditional probability twice to isolate the event $t+1$ is pivotal into the component pieces $\tau_m = t+1$ for each $m$,
    \begin{equs}\label{eq:Baye-two-times}
        ~&\pr(C_{t+1}\geq \frac{2}{3}L_0|t+1 \text{ is pivotal, }T_i\leq t < T_{i+1})\\
        &= \frac{\pr(C_{t+1}\geq \frac{2}{3}L_0, t+1 \text{ is pivotal}|T_i\leq t < T_{i+1})}{\pr(t+1 \text{ is pivotal}|T_i\leq t < T_{i+1})}\\
        &=\frac{1}{\pr(t+1 \text{ is pivotal}|T_i\leq t < T_{i+1})}\sum_{m \leq n}  \pr(C_{t+1}\geq \frac{2}{3}L_0,\, t+1 = \tau_m| T_i\leq t < T_{i+1})\\
        &=\frac{1}{\pr(t+1 \text{ is pivotal}|T_i\leq t < T_{i+1})}\sum_{m \leq n}  \pr(C_{\tau_m}\geq \frac{2}{3}L_0| \tau_m=t+1,\, T_i\leq t < T_{i+1})\pr(\tau_m=t+1|T_i\leq t < T_{i+1}).
    \end{equs}
    Now observe that $\{\tau_m=t+1 , T_i\leq t < T_{i+1}\} \in \Lambda_{\tau_m}$, so applying the Markov property for $(\Lambda_t)$, Lemma \ref{lm:New Level Markov Property},
    \begin{equs}\label{eq:big-cycle-prob}
        ~&\pr(C_{\tau_m}\geq \frac{2}{3}L_0| \tau_m=t+1 , T_i\leq t < T_{i+1})\\
        &=\frac{\E[\pr_{\pi \sim \mrm{Uniform}(S_{\lambda_{\tau_m}})}(n_{\tau_m} \text{ is in a cycle of length } \geq \frac{2}{3}L_0)\mathbbm{1}_{\tau_m=t+1 , T_i\leq t < T_{i+1}}]}{\E[\mathbbm{1}_{\tau_m=t+1 , T_i\leq t < T_{i+1}}]}\\
        &\geq \frac{2(k-i)L_0}{3n},
    \end{equs}
    as on the event $\{\tau_m=t+1 , T_i\leq t < T_{i+1}\}$, $\lambda_{\tau_m}$ has $k-i$ parts of size at least $\frac{2}{3}L_0$, and the probability $n_{\tau_m}$ belongs to a cycle corresponding to any part $(\lambda_{\tau_m})_i$, is exactly $\frac{(\lambda_{\tau_m})_i}{n}$ under the uniform distribution of $S_{\lambda_{\tau_m}}$. So finally plugging \eqref{eq:big-cycle-prob} back into \eqref{eq:Baye-two-times}, and using the fact that,
    \begin{equs}
        \sum_{m }\pr(\tau_m=t+1|T_i\leq t < T_{i+1})=\pr(t+1 \text{ is pivotal}|T_i\leq t < T_{i+1}),
    \end{equs}
    we recover \eqref{eq:See-big-cycle-at-pivot}.
    %The first step to seeing this is to observe that,
    %\begin{equs}
    %    \pr(T^{(i)}>t| \sigma_0,\sigma_1,\dots,\sigma_{n+|\sigma|}) = \frac{\#\{\text{cycles of } \sigma_t \text{ of length }\geq \frac{2}{3}L_0\}}{k}.
    %\end{equs}
    %This is because $\mrm{c}_t^{(i)}$ is distributed uniformly among the cycles of $\sigma_t$ of length $\geq \frac{2}{3}L_0$ \todo{elaborate}. So using Baye's theorem \todo{FIXXXX}
    %\begin{equs}
    %    ~&\pr(\mf{c}_{t+1}^{(i)} \text{is the pivotal cycle at time } t+1| \{t+1 \text{ is a pivotal time}\} \cap \{T^{(i)}>t\})\\
    %    &= \pr(\mf{c}_{t+1}^{(i)} \text{is the pivotal cycle at time } t+1|t+1 \text{ is a pivotal time})\\
     %   &\cdot\frac{\pr(T^{(i)}>t| \{t+1 \text{ is a pivotal time}\} \cap \{\mf{c}_{t+1}^{(i)} \text{is the pivotal cycle at time } t+1\})}{\pr(T^{(i)}>t|t+1 \text{ is a pivotal time})}\\
    %    &=\pr(\mf{c}_{t+1}^{(i)} \text{is the pivotal cycle at time } t+1|t+1 \text{ is a pivotal time})
    %\end{equs}

    %\begin{equs}
     %   =\frac{\mrm{length}(\mf{c}_{t+1}^{(i)})}{n} \geq \frac{2L_0}{3n}
    %\end{equs}
    %\todo{finish}.
    \textbf{Step 3:} Next, we need to prove that,
    \begin{equs}\label{eq:pivotal-cycle-splits}
        \pr(t_{i+1} \leq t+3| T_i \leq t < T_{i+1}) \geq \frac{(k-i)L_0^{1/2}}{240n}.
    \end{equs}
    As with the conditional probability arguments of the last step, it suffices to show that,
    \begin{equs}\label{eq:pivotal-cycle-splits-2}
        ~&\pr(t_{i+1} \leq t+2| C_{t+1} \geq \frac{2}{3}L_0, T_i \leq t < T_{i+1}) \geq \frac{1}{40L_0^{1/2}}, \text{ and}\\
        &\pr(t_{i+1} \leq t+3| C_{t+2} \geq \frac{2}{3}L_0, T_i \leq t < T_{i+1}) \geq \frac{1}{40L_0^{1/2}}.
    \end{equs}
    Note $\frac{1}{40L_0^{1/2}} \times \frac{(k-i)L_0}{6n}=\frac{(k-i)L_0^{1/2}}{240n}$. And again in parallel with the last step, we will only prove the first equation of \eqref{eq:pivotal-cycle-splits-2}, as the second is proven identically. Now observe that if the pivotal cycle at time $t+1$ has length $\geq \frac{2}{3}L_0$, and in the next step it splits into two cycles, which each have length less than $2/3$ the original size, then the number of cycles of size $\geq \frac{2}{3}L_0$ has reduced by $1$ from $\sigma_{t+1}$ to $\sigma_{t+2}$. In other words, 
    \begin{equs}
        \pr(t_{i+1} \leq t+2| C_{t+1} \geq \frac{2}{3}L_0, T_i \leq t < T_{i+1})\geq \pr(C_{t+2} \in (\frac{1}{3}C_{t+1},\frac{2}{3}C_{t+1}) | C_{t+1} \geq \frac{2}{3}L_0, T_i \leq t < T_{i+1}).
    \end{equs}
    But by the Markov property and definition of the Weingarten process with respect to $\mc{F}_t$,
    \begin{equs}
        \pr(C_{t+2} \in (\frac{1}{3}C_{t+1},\frac{2}{3}C_{t+1}) | \mc{F}_t) &= \sum_{j \in (\frac{1}{3}C_{t+1},\frac{2}{3}C_{t+1})} \frac{\mrm{Cat}(j-1)\mrm{Cat}(C_{t+1}-j-1)}{\mrm{Cat}(j-1)}\\
        &\geq \sum_{j \in (\frac{1}{3}C_{t+1},\frac{2}{3}C_{t+1})} \frac{1}{6C_{t+1}^{3/2}} \geq \frac{1}{24 C_{t+1}^{1/2}}\mathbbm{1}_{C_{t+1} \geq 4},
    \end{equs}
    where we used the Catalan number quotient bound, Lemma \ref{Catalan quotient bound}, in the second to last inequality, and in the last inequality, we assumed $C_{t+1} \geq 4$ so that at least $1/4$ of the values in $[C_{t+1}]$ lie in the open interval $(\frac{1}{3}C_{t+1},\frac{2}{3}C_{t+1})$. Thus on the event, $\{C_{t+1} \geq \frac{2}{3}L_0, T_i \leq t < T_{i+1}\}$, $ \pr(C_{t+2} \in (\frac{1}{3}C_{t+1},\frac{2}{3}C_{t+1}) | \mc{F}_t) \geq \frac{1}{36L_{0}^{1/2}}\geq \frac{1}{40L_{0}^{1/2}}$, and so averaging on this event completes the proof of \eqref{eq:pivotal-cycle-splits-2} (recall we assumed $L_0 \geq 6$, so $C_{t+1}\geq 4$).
    %Condition on $\sigma_{t+1}^{(i)}=\pi$ satisfying $\mf{c}_{t+1}^{(i)} \text{is the pivotal cycle at time } t+1$ and the length of this pivotal cycle in $\pi$ is $\ell \geq \frac{2L_0}{3n}$. Then,
    %\begin{equs}
    %    \pr(t^{(i)} = t+2| \sigma_{t+1}^{(i)}=\pi)&\leq \pr(C_{t+2} \in (\frac{1}{3}\ell,\frac{2}{3}\ell)| \sigma_{t+1}^{(i)}=\pi)\\
     %   &= \sum_{j \in (\frac{1}{3}\ell,\frac{2}{3}\ell)} \frac{\mrm{Cat}(j-1)\mrm{Cat}(\ell-j-1)}{\mrm{Cat}(\ell-1)}\\
     %   &\geq \sum_{j \in (\frac{1}{3}\ell,\frac{2}{3}\ell)} \frac{1}{6\ell^{3/2}}\\
     %   &\geq \frac{1}{24\ell^{1/2}} \geq \frac{1}{40L_0^{1/2}},
    %\end{equs}
    %where we have,\todo{Explain!} Similarly conditioning on $\sigma_{t+2}^{(i)}=\pi$ satisfying $\mf{c}_{t+2}^{(i)} \text{is the pivotal cycle at time } t+2$,
    %\begin{equs}
     %   \pr(t^{(i)} = t+3| \sigma_{t+2}^{(i)}=\pi)\geq \frac{1}{40L_0^{1/2}}.
    %\end{equs}
    %It follows that
    %\begin{equs}
     %   \pr(t^{(i)} \leq t+3| \{T^{(i)}> t\} \cap \mf{c}_{t+1}^{(i)} \text{ or } \mf{c}_{t+2}^{(i)} \text{ are pivotal cycles}\})\geq \frac{1}{40L_0^{1/2}},
    %\end{equs}
    %and hence combining with \eqref{eq:pivotal-cycle-appears}, we obtain \eqref{eq:pivotal-cycle-splits}. 
    
    \textbf{Step 4:} Finally as with \eqref{eq:Pivotal in two steps}, we have
    \begin{equs}
        \pr(t+4 \text{ or } t+5 \text{ are pivotal times}| \mc{F}_{t+3}) \geq \frac{1}{4},
    \end{equs}
    and so at once we have,
    \begin{equs}
        \pr(T_{i+1} \leq t+5| t_i \leq t+3)\geq \frac{1}{4},
    \end{equs}
    and as a result,
    \begin{equs}
        \pr(T_{i+1} \leq t+5| T_i \leq t <T_{i+1}) \geq \frac{1}{4} \times \frac{(k-i)L_0^{1/2}}{240n} \geq 10^{-3}\frac{(k-i)L_0^{1/2}}{n}.
    \end{equs}
    By direct corollary, we see inductively that
    \begin{equs}
        \pr(T_{i+1}-T_i \geq 5t) &= \pr(T_{i+1}-T_i \geq 5t | T_{i+1}-T_i \geq 5(t-1))\cdots \pr(T_{i+1}-T_i \geq 5 | T_{i+1}-T_i \geq 0)\\
        &\leq (1-10^{-3}(k-i)L_0^{1/2}n^{-1})^t.
    \end{equs}
\end{proof}

We are now ready to prove Lemma \ref{lm:time-to-halve}. The key is to treat all of the long cycles of length $\geq \frac{2}{3}L_0$ as coupons and to bound $T$ similarly to the classical coupon collector problem.

\begin{proof}[Proof of Lemma \ref{lm:time-to-halve}]
    Using a union bound, \eqref{eq:T^i tail}, and the fact that $1+\log k \geq 1+\frac{1}{2}+\dots+\frac{1}{k}$,
    \begin{equs}
        \pr(T>5t(1+\log k))\leq \pr\bigg(\bigcup_{i=0}^{k-1}\{T_{i+1}-T_i \geq 5\frac{t}{k-i}\} \bigg)
        &\leq \sum_{i=0}^{k-1} \bigg(1-10^{-3}\frac{(k-i)L_0^{1/2}}{n}\bigg)^{\frac{t}{k-i}}.
    \end{equs}
    Now substituting $t \mapsto 10^3 \frac{n}{L_0^{1/2}}(t+\log k)$, applying the elementary inequality $(1-x^{-1})^x \leq e^{-1}$ for $x >1$, ,
    \begin{equs}
        \pr\bigg(T>5000\frac{n}{L_0^{1/2}}(1+\log k)(t+\log k)\bigg) \leq \sum_{i=0}^{k-1} e^{-t-\log k} = e^{-t}.
    \end{equs}
    Finally recalling that the number of comparatively long cycles at time zero, $k$, is bounded by $\frac{n}{\frac{2}{3}L_0}$ we have,
    \begin{equs}\label{eq:T-tail-bound}
        \pr\bigg(T>5000\frac{n}{L_0^{1/2}}\bigg(1+\log\bigg(\frac{3n}{2L_0}\bigg)\bigg)\bigg(t+\log\bigg(\frac{3n}{2L_0}\bigg)\bigg)\leq e^{-t}.
    \end{equs}
    So decomposing $T$ via the cutoff scale $T_{\mrm{cutoff}} :=5000\frac{n}{L_0^{1/2}}\bigg(1+\log\bigg(\frac{3n}{2L_0}\bigg)\bigg)\log\bigg(\frac{3n}{2L_0}\bigg)$, and applying the layer cake formula for expectations,
    \begin{equs}
        \E[T] &= \E\bigg[T \mathbbm{1}\{T \leq T_{\mrm{cutoff}}\}\bigg]+\E\bigg[T \mathbbm{1}\{T > T_{\mrm{cutoff}}\}\bigg]\\
        &\leq T_{\mrm{cutoff}}+\int_{0}^{\infty}\pr(T>T_{\mrm{cutoff}}+u)du\\
        &\leq T_{\mrm{cutoff}}+5000\frac{n}{L_0^{1/2}}\bigg(1+\log\bigg(\frac{3n}{2L_0}\bigg)\bigg)\int_{0}^{\infty} \pr\bigg(T>5000\frac{n}{L_0^{1/2}}\bigg(1+\log\bigg(\frac{3n}{2L_0}\bigg)\bigg)\bigg(t+\log\bigg(\frac{3n}{2L_0}\bigg)\bigg) dt\\
        &\leq T_{\mrm{cutoff}}+5000\frac{n}{L_0^{1/2}}\bigg(1+\log\bigg(\frac{3n}{2L_0}\bigg)\bigg)\int_{0}^{\infty} e^{-t} dt\\
        &\leq 10^4\frac{n}{L_0^{1/2}}\bigg(1+\log\bigg(\frac{3n}{2L_0}\bigg)\bigg)^2,
    \end{equs}
    where we used the change of variables $t=5000\frac{n}{L_0^{1/2}}\bigg(1+\log\bigg(\frac{3n}{2L_0}\bigg)\bigg)u$, and the bound \eqref{eq:T-tail-bound} in the equation display above.
\end{proof}
\color{black}
\section{Unitary Group: Additional Estimates}\label{Small Permutations Section}

This section relies on definitions from subsections \ref{U Weingarten Graph Subsection} and \ref{U Weingarten Function Section}.

\subsection{Outline of the Proofs of Theorems \ref{Small Permutations Theorem} and \ref{n^{5/4} Theorem}}

The proof of Theorem \ref{Small Permutations Theorem} will be a combinatorial analysis of $|\mathrm{P}(\sigma,|\sigma|+2g)|$ using elementary tools such as the ``stars and bars'' method. The analysis starts by separating the solid edges of a path $\mathrm{p} \in \mathrm{P}(\sigma,|\sigma|+2g)$ into "defect" and "non-defect" edges, where non-defect edges are the edges taking a step in the direction of the identity. The key insight is that once we select the ordering of the non-defect solid edges, the placement and choice of transpositions associated to the defect edges are severely restricted for permutations with $|\sigma|$ small. The bound on the Weingarten function will the follow from plugging into the path expansion formula.

The proof of Theorem \ref{n^{5/4} Theorem} will use Theorem \ref{Small Permutations Theorem}. The approach is based on the  ``loop equation'' satisfied by the Weingarten function (equivalent to the $\mathrm{U}(N)$ lattice gauge theory master loop equation \cite{AN2023,CPS2023,borga2024surface}). This equation can be written in an operator theoretic way on a vector space with components indexed by permutations. Then with the correct norm on this vector space, we will obtain a "forcing" estimate similar to energy estimates coming from PDE theory (see \cite{borga2024surface,CNS2025a} for a similar argument). Through this forcing estimate which uses the small permutation regime result, we will be able to prove Theorem \ref{n^{5/4} Theorem}.

\subsection{Small Permutation Regime for Weingarten Function}\label{Small Permutations}

\begin{lemma}\label{Small Permutation Path Count Bound}
    There is a constant $C_{|\sigma|}$ only depending on $|\sigma|$, and in particular not on $n$, such that for all $\sigma \in S_n$ and $g \geq 0$,
    \begin{equs}
        |\mathrm{P}(\sigma,|\sigma|+2g)| \leq C_{|\sigma|}(\sqrt{48}e n)^{2g}.
    \end{equs}
\end{lemma}

\begin{proof}
    Fix any path $\mathrm{p}=(\sigma_1,\sigma_2,\dots\sigma_{n+|\sigma|+2g}) \in \mathrm{P}(\sigma,|\sigma|+2g)$. We call a solid edge in the path \textit{defective} if it takes $\sigma_i \to \sigma_{i+1}$ with $|\sigma_{i+1}|>|\sigma_i|$. We can observe that $\mathrm{p}$ must have exactly $g$ defective solid edges, and $|\sigma|+g$ non-defective solid edges.  

    Any path $\mathrm{p}\in \mathrm{P}(\sigma,|\sigma|+2g)$ can then be constructed, by first abstractly, taking a sequence of $n$ dashed edges, $g$ defect solid edges, and $g+|\sigma|$ non-defect solid edges some order, and assigning corresponding transpositions to each solid edge. To count the number of valid constructions as detailed above, start with a sequence of $n$ dashed edges and then choose where in the sequence to insert the $g$ defective solid edges, there are at most $\binom{n+g}{g}$ choices for this placement. Next we need to choose the transposition assigned to each of these defective solid edges. Since we already chose the levels where these transpositions occur, there are at most $n$ choices for each transposition, and therefore there are at most $n^g$ choices for the $g$ defective transpositions. 

    Now we still have to choose the placement of the non-defective solid edges. The key observation is that if $g+|\sigma|<n/2$ with only $g$ defect solid edges, there must be an $n-2g-2|\sigma|$ element subset $X$ of $[n]$, such that any permutation in $\mathrm{p}$ acts as the identity on $X$. This is because any permutation $\sigma$ has at most $n-2|\sigma|$ fixed points, and any of the $g$ defect solid edge step, can add at most 2 more elements which are not fixed by every permutation in the path, while any other step of the path, cannot change to this number. Thus, there are at most $2g+2|\sigma|$ levels where the non-defective solid edges can be placed, and all together these levels have at most $g$ solid edges occurring in $\mathrm{p}$ at these levels. As a result, there are at most $\binom{(2g+2|\sigma|)+g+(g+|\sigma)|}{g+|\sigma|}$ choices for the placement of the non-defective solid edges. Finally as all the permutations in the path act as the identity on $X$, there are at most $(2g+2|\sigma|)^{|\sigma|+g}$ choices for the transpositions assigned to the non-defective solid edges. Putting it all together, if $g+|\sigma| <n/2$, then
    \begin{equs}
        |\mathrm{P}(\sigma,|\sigma|+2g)| &\leq \binom{n+g}{g}n^g \binom{4g+3|\sigma|}{g+|\sigma|} (2g+2|\sigma|)^{g+|\sigma|}\\
        &\leq \frac{(3/2)^g n^g}{g!} n^g 2^{4g+3|\sigma|} 2^{g+|\sigma|} (g+|\sigma|)^{g+|\sigma|}   \\
        &= 48^g 2^{4|\sigma|} \frac{g^g}{g!}g^{|\sigma|}(1+|\sigma|/g)^{g+|\sigma|} n^{2g}\\
        &\leq 48^g 2^{4|\sigma|} e^{g}(|\sigma|)!e^g e^{|\sigma|+|\sigma|^2/g} n^{2g}\\
        &\leq e^{|\sigma|^2}(16 e g)^{|\sigma|}(|\sigma|)!(\sqrt{48}e n)^{2g},
    \end{equs}
    where we used the binomial theorem to bound $\binom{4g+3|\sigma|}{g+|\sigma|} \leq 2^{4g+3|\sigma|}$, the Stirling bound $g^g \leq e^g g!$ and $g^{|\sigma|}\leq (|\sigma|)!e^g$, and the bound $1+|\sigma|/g \leq e^{|\sigma|/g}$.

    On the other hand, if $g+|\sigma| \geq n/2$, then there are are still at most $\binom{4g+3|\sigma|}{g+|\sigma|}$ choices for the placement of the non-defective solid edges, but the number of choices for the transpositions they are assigned can be bounded by $n^g (2g+2|\sigma|)^{|\sigma|}$ instead, so we have,
    \begin{equs}
        |\mathrm{P}(\sigma,|\sigma|+2g)| &\leq \binom{n+g}{g}n^g \binom{4g+3|\sigma|}{g+|\sigma|} n^g (2g+2|\sigma|)^{|\sigma|}\\
        &\leq 2^{3g/2} n^g 2^{4g+3|\sigma|} n^g (2 g)^{|\sigma|}(1+|\sigma|/g)^{|\sigma|}\\
        &\leq 2^{11g/2}(16g)^{|\sigma|} e^{|\sigma|^2}n^{2g}.
    \end{equs}
    for $g+|\sigma| \geq n/2$. Finally $g^{|\sigma|} \leq (|\sigma|)! e^{g}$ so $|\mathrm{P}(\sigma,|\sigma|+2g)| \leq (|\sigma|)!16^{|\sigma|}e^{|\sigma|^2}(2^{11/2}e)^g n^{2g}$ for $g+|\sigma| \geq n/2$.

    Combining the two regimes we get the desired bound.
\end{proof}

\begin{proof}[Proof of Theorem \ref{Small Permutations Theorem}]
    Applying the previous lemma, for $N > \sqrt{48}en$
    \begin{equs}
        |\overline{\Wg}_N(\sigma)-\mathrm{Moeb}(\sigma)| &= \sum_{g \geq 1} |\mathrm{P}(\sigma,|\sigma|+2g)|N^{-2g}\\
        &\leq C_{|\sigma|}\sum_{g \geq 1} (\sqrt{48}en/N)^{2g}\\
        &=\frac{C_{|\sigma|}(48e^2n^2/N^2)}{1-48e^2n^2/N^2}\\
        &= \frac{\tilde{C}_{|\sigma|}n^2}{N^2-48e^2n^2},
    \end{equs}
    for $\tilde{C}_{|\sigma|}=48e^2C_{|\sigma|}$.
\end{proof}

\subsection{Energy estimate for the Loop Equation}\label{Fixed Point Section}
For this section let $\mathcal{W}_n$ denote the unitary Weingarten graph up to level $n$ so that $\mathcal{W}_n$ has vertex set $\bigsqcup_{0\leq k \leq n} S_k$. Next let $\R^{\mathcal{W}_n}$ denote the $1+\sum_{k=1}^{n}k!$ dimensional vector space with the elements of $\mathcal{W}_n$ as basis vectors. We can reformulate the orthogonality recursion for Weingarten functions derived in \cite{Collins2017} as the following fixed point problem in $\R^{\mathcal{W}_n}$.

\begin{lemma}
    Suppose $N \geq n$. Then the vector $w := (\Wg_N(\sigma))_{\sigma \in \mathcal{W}_n}$ solves the fixed point problem,
    \begin{equs}
        w=Tw,
    \end{equs}
    where $T: \R^{\mathcal{W}_n} \to \R^{\mathcal{W}_n}$ is the linear operator given by $(Tx)_{\emptyset} = x_{\emptyset}$, and,
    \begin{equs}
        (Tx)_{\sigma}=N^{-1}x_{\sigma^{\downarrow}} \delta_{\sigma(k),k}-N^{-1}\sum_{i=1}^{k-1} x_{(i \, k)\sigma},
    \end{equs}
    for each $\sigma \in S_k \subset \mathcal{W}_n$ ($k \neq 0$), where $\sigma^{\downarrow}$ is $\sigma$ restricted to a permutation in $S_{k-1}$.
\end{lemma}

Our goal will be to use the fixed point formulation above, to get quantitative bounds on $w$. To preform our analysis, we will need to first reformulate our fixed point problem in the space of conjugacy classes of permutations. 

 Recalling that the Weingarten function is a class function, we can reformulate our fixed point problem on $\R^{\mathfrak{P}_n}$ where $\mathfrak{P}_n$ is the set of partitions $\lambda \vdash r$ for $r \leq n$. We will sometimes abuse notation by referring to the component $x_{\lambda}$ of a vector $(x_{\lambda}) \in \R^{\mathfrak{P}_n}$ by $x_{\sigma}$ for some $\sigma$ in the conjugacy class $S_{\lambda}$. Lastly for any conjugacy class of $S_k$, $S_\lambda$, we will assign a permutation $\sigma_{\lambda} \in S_{\lambda}$ such that $k$ is in the largest cycle of $\sigma_{\lambda}$.

\begin{lemma}[Reformulated Fixed point problem]
    Suppose $N \geq n$. Then the vector $w := (\Wg_N(\lambda))_{\lambda\in \mathfrak{P}_n}$ solves the  fixed point problem,
    \begin{equs}
        w=\Tilde{T}w,
    \end{equs}
    where $\Tilde{T}: \R^{\mathfrak{P}_n} \to \R^{\mathfrak{P}_n}$ is the operator given by $(\Tilde{T}x)_{\emptyset} = x_{\emptyset}$, and,
    \begin{equs}\label{eq:Partition-fixed-point-problem}
        (\Tilde{T}x)_{\lambda} = N^{-1}\Wg_N(\groupid_{k-1})\delta_{\sigma_{\lambda},\groupid_k}-N^{-1}\sum_{i=1}^{k-1} x_{(i \, k)\sigma_{\lambda}},
    \end{equs}
    where $\groupid_k$ denotes the identity permutation in $S_k$.
\end{lemma}
\begin{proof}
    Since $\sigma_{\lambda} \in S_{\lambda} \in S_k$ is a permutation with the largest cycle containing $k$. $\sigma_{\lambda}$, fixes $k$ if and only if $\sigma_{\lambda}=\groupid_{k}$. Therefore the Weingarten recursion from \cite{Collins2017} exactly reads
    \begin{equs}              \Wg_N(\groupid_{r})=N^{-1}\Wg_N(\groupid_{r-1})\delta_{\sigma_{\lambda},\groupid_r}-N^{-1}\sum_{i=1}^{r-1} \Wg_N((i \, k)\sigma_{\lambda}),
    \end{equs}
    thus $(\Wg_N(\lambda))_{\lambda\in \mathfrak{P}_n}$ satisfies \eqref{eq:Partition-fixed-point-problem}.
\end{proof}
We now define a norm on $\R^{\mathfrak{P}_n}$.
\begin{definition}
    For $x=(x_{\lambda})_{\lambda \in \mathfrak{P}_n}$, and $\gamma>0$, we define the norm of $x$ as follows,
    \begin{equs}
        \|x\|_{\gamma}:= |x_{\emptyset}|+ \sup_{r \in [n]} \sup_{\lambda \vdash r} N^{r+|\sigma_{\lambda}|} \gamma^{|\sigma_{\lambda}|} |\mathrm{Moeb}(\sigma_{\lambda})|^{-1}|x_{\lambda}|.
    \end{equs}
\end{definition}

With this setup in place, we now prove the key estimate of the section.

\begin{lemma}\label{PDE Type Inequality}
    For $n<\frac{1}{100\sqrt{48}e}N$, and for any $x \in \R^{\mathfrak{P}_n}$, we have the bound
    \begin{equs}
        \|\tilde{T}x\|_{\gamma} \leq \bigg(\gamma+\frac{12n^{5/2}}{N^2}\bigg)\|x\|_{\gamma}+1+\frac{n^2}{50N^2},
    \end{equs}
    for $\gamma \in (\frac{1}{2},1)$.
\end{lemma}

\begin{proof}
    Recall for each $\lambda \vdash r$ for $r \geq 1$,
    \begin{equs}
        (\tilde{T}x)_{\lambda}:= N^{-1}\Wg_N(\groupid_{r-1})\delta_{\sigma_{\lambda},\groupid_r}-N^{-1}\sum_{i=1}^{r-1} x_{(i \, r)\sigma_{\lambda}}, 
    \end{equs}
    and as a result, denoting $\sigma$ by $\sigma_{\lambda}$ for notational simplicity, we have,
    \begin{equs}\label{eq:fixed-point-1}
        &N^{r+|\sigma|} \gamma^{|\sigma|} |\mathrm{Moeb}(\sigma)|^{-1} |(\tilde{T}x)_{\lambda}|\\
        &\leq N^{r-1+|\sigma|} \gamma^{|\sigma|} |\mathrm{Moeb}(\sigma)|^{-1}|\Wg_N(\groupid_{r-1})|\delta_{\sigma,\groupid_r}+N^{r-1+|\sigma|} \gamma^{|\sigma|} |\mathrm{Moeb}(\sigma)|^{-1}\sum_{i=1}^{r-1} x_{(i \, r)\sigma}.
    \end{equs}
    But $N^{r-1+|\sigma|} \gamma^{|\sigma|} |\mathrm{Moeb}(\sigma)|^{-1}|\Wg_N(\groupid_{r-1})|\delta_{\sigma,\groupid_r}=0$ unless $\sigma=\groupid$ in which case it is equal to $N^{r-1} |\Wg_N(\groupid_{r-1})|$. Additionally, we can split the sum $\sum_{i=1}^{r-1} x_{(i \, r)\sigma}$ into the terms in the sum which split a cycle of $\sigma$, and those which merge two cycles. Lastly observe that in the splitting terms, $|(i \, r)\sigma|=|\sigma|-1$, while for the merger terms $|(i \, r)\sigma|=|\sigma|+1$. These observations together with \eqref{eq:fixed-point-1}, and rewriting $\mathrm{Moeb}(\sigma)=\frac{\mathrm{Moeb}(\sigma)}{\mathrm{Moeb}((i \, r)\sigma)}\mathrm{Moeb}((i \, r)\sigma)$, allow us to bound,
    \begin{equs}
        N^{r+|\sigma|} \gamma^{|\sigma|} |\mathrm{Moeb}(\sigma)|^{-1} |(\tilde{T}x)_{\lambda}| \leq a_1+a_2+a_3,
    \end{equs}
    with
    \begin{equs}
        ~&a_1:= N^{r-1} |\Wg_N(\groupid_{r-1})|,\\
        &a_2:= \gamma \sum_{i \sim_{\sigma} r} \frac{|\mathrm{Moeb}((i \, r)\sigma)|}{|\mathrm{Moeb}(\sigma)|} N^{r+|(i \, r)\sigma|} \gamma^{|(i \, r)\sigma|}|\mathrm{Moeb}((i \, r)\sigma)|^{-1}|x_{(i \, r)\sigma}|,\\
        &a_3:= N^{-2}\gamma^{-1}\sum_{i \nsim_{\sigma} r} \frac{|\mathrm{Moeb}((i \, r)\sigma)|}{|\mathrm{Moeb}(\sigma)|} N^{r+|(i \, r)\sigma|} \gamma^{|(i \, r)\sigma|}|\mathrm{Moeb}((i \, r)\sigma)|^{-1}|x_{(i \, r)\sigma}|,
    \end{equs}
    where $i \sim_{\sigma} r$ ($i \sim_{\sigma} r$) means that $i$ and $r$ (do not) belong to the same cycle in $\sigma$.
    Now observe as a direct consequence of Theorem \ref{Small Permutations Theorem}, that
    \begin{equs}
        &a_1= N^{r-1} |\Wg_N(\groupid_{r-1})| \leq 1+\frac{n^2}{50N^2}.
    \end{equs}
    Next, notice that $\sum_{i \sim_{\sigma} r}|\mathrm{Moeb}((i \, r)\sigma)|=|\mathrm{Moeb}(\sigma)|$. This is a consequence of the Catalan number recursion since the product of Catalan numbers resulting from the cycles in $\sigma$ not containing $r$ are present on both size, while the sum is over exactly all ways to split the cycle containing $r$ into two parts. We use the definition of the norm $\|\cdot\|_{\gamma}$ together with the observation of the previous sentence to bound $a_2$,
    \begin{equs}
        a_2 &= \gamma \sum_{i \sim_{\sigma} r} \frac{|\mathrm{Moeb}((i \, r)\sigma)|}{|\mathrm{Moeb}(\sigma)|} N^{r+|(i \, r)\sigma|} \gamma^{|(i \, r)\sigma|}|\mathrm{Moeb}((i \, r)\sigma)|^{-1}|x_{(i \, r)\sigma}|\\
        &\leq \gamma \sum_{i \sim_{\sigma} r} \frac{|\mathrm{Moeb}((i \, r)\sigma)|}{|\mathrm{Moeb}(\sigma)|} \|x\|_{\gamma}\\
        &= \gamma \|x\|_{\gamma}.
    \end{equs}

    Finally, to bound $a_3$, we recall as in the proof of Lemma \ref{Single Defect Lemma}, that for $i \nsim_{\sigma_{\lambda}} r$, $\frac{|\mathrm{Moeb}((i \, r)\sigma)|}{|\mathrm{Moeb}(\sigma)|}$ is exactly a quotient of Catalan numbers where Lemma \ref{Catalan quotient bound} directly applies and we have,
    \begin{equs}
        a_3 &= N^{-2}\gamma^{-1}\sum_{i \nsim_{\sigma} r} \frac{|\mathrm{Moeb}((i \, r)\sigma)|}{|\mathrm{Moeb}(\sigma)|} N^{r+|(i \, r)\sigma|} \gamma^{|(i \, r)\sigma|}|\mathrm{Moeb}((i \, r)\sigma)|^{-1}|x_{(i \, r)\sigma_{\lambda}}|\\
        & \leq N^{-2}\gamma^{-1} \sum_{i \nsim_{\sigma} r} \frac{|\mathrm{Moeb}((i \, r)\sigma)|}{|\mathrm{Moeb}(\sigma)|}\|x\|_{\gamma}\\
        &\leq n N^{-2}\gamma^{-1} \sup_{i \nsim_{\sigma} r} \frac{|\mathrm{Moeb}((i \, r)\sigma)|}{|\mathrm{Moeb}(\sigma)|}\|x\|_{\gamma}\\
        &\leq 6 n^{5/2} N^{-2}\gamma^{-1} \|x\|_{\gamma}.
    \end{equs}
    Combining all the bounds obtained so far and taking the supremum over all $\lambda \in \mf{P}_n$, we have,
    \begin{equs}
        \|\tilde{T}x\|_{\gamma} \leq \bigg(\gamma+\frac{6\gamma^{-1}n^{5/2}}{N^2}\bigg)\|x\|_{\gamma}+1+\frac{n^2}{50N^2},
    \end{equs}
    and now restricting to $\gamma \in (1/2,1)$, we obtain the desired bound.
\end{proof}

We are now ready to prove the main theorem of the section.

\begin{theorem}
    Suppose $n \leq 10^{-4}N^{4/5}$. Then for any $\sigma \in S_n$ with $|\sigma| \geq 4$,
    \begin{equs}\label{eq:Fixed-point-conclusion}
        \frac{\overline{\Wg}_N(\sigma)}{\mathrm{Moeb}(\sigma)}\leq  e^2|\sigma| \exp \bigg(25|\sigma|\frac{n^{5/2}}{N^2}\bigg),
    \end{equs}
\end{theorem}

\begin{proof}
    Let $w=(\Wg_N(\lambda))_{\lambda\in \mathfrak{P}_n}$ denote the fixed point of \eqref{eq:Partition-fixed-point-problem}, and directly applying Lemma \ref{PDE Type Inequality} to the fixed point $x=w$, we see that for any $\gamma \in (1/2,1)$,
    \begin{equs}
        ||w||_{\gamma} \leq \bigg(\gamma+\frac{12n^{5/2}}{N^2}\bigg)||w||_{\gamma}+1+\frac{n^2}{50N^2},
    \end{equs}
    so after setting $\gamma=1-\frac{12n^{5/2}}{N^2}-|\sigma|^{-1}$ and rearranging the inequality,
    \begin{equs}
        ||w||_{1-\frac{12n^{5/2}}{N^2}-|\sigma|^{-1}} \leq |\sigma|\bigg(1+\frac{n^2}{50N^2}\bigg).
    \end{equs}
    Plugging in the definition of $||w||_{1-\frac{12n^{5/2}}{N^2}-|\sigma|^{-1}}$,
    \begin{equs}
        \frac{N^{n+|\sigma|}|\mathrm{Wg}_N(\sigma)|}{|\mathrm{Moeb}(\sigma)|} &\leq \bigg(1-\frac{12n^{5/2}}{N^2}-|\sigma|^{-1}\bigg)^{-|\sigma|}|\sigma|\bigg(1+\frac{n^2}{50N^2}\bigg)\\
        &\leq \exp\bigg(2\bigg(\frac{12n^{5/2}}{N^2}+|\sigma|^{-1}\bigg)|\sigma|\bigg)|\sigma|\exp\bigg(\frac{n^2}{50N^2}\bigg)\\
        &\leq e^2|\sigma| \exp \bigg(25|\sigma|\frac{n^{5/2}}{N^2}\bigg),
    \end{equs}
    where the second to last line results from the elementary inequalities, $1+x \leq e^x$ for $x \in \R$, and $(1-x)^{-1}\leq e^{2x}$, for $x \in (0,1/2)$.
\end{proof}
Theorem \ref{n^{5/4} Theorem} now follows as a direct corollary by taking the logarithm in equation \eqref{eq:Fixed-point-conclusion}.

\section{Orthogonal Groups: Proof of Theorem \ref{Main Bound Orthogonal}}\label{Orthogonal Section}

We let $\mc{P}_2(2n)$ denote the set of pairings of the set $2n$, that is,
\begin{equs}
    \mc{P}_2(2n) = \{\{\{k_1,k_2\},\{k_3,k_4\},\dots,\{k_{2n-1},k_{2n}\}\}: \{k_1,k_2,\dots,k_{2n}\} = [2n]\}.
\end{equs}
For $G \in \{\mathrm{O}(N),\mathrm{SP}(N)\}$, the function $\Wg_N^{G} (\pi)$ takes in pairings $\pi \in \mathcal{P}_2(2n)$ as inputs. The main result of this section is Theorem \ref{Main Bound Orthogonal} which we restate here: For $G \in \{\mathrm{O}(N),\mathrm{SP}(N)\}$, there is a constant $C>0$ such that,
    \begin{equs}
        \frac{N^{-n-|\pi|}|\mathrm{Wg}_N^{G}(\pi)|}{\mathrm{Moeb}(\pi)}\leq \frac{1}{(1-Cn^{3/2}/N)^2}.
    \end{equs}
    For any $\pi \in \mathcal{P}_2(2n)$ and $Cn^3<N^2$.

The definitions of $|\pi|$ and $\mathrm{Moeb}(\pi)$ will be given in the next subsection.

\begin{remark}\label{Orthogonal to Symplectic}
    We will restrict the remainder of this section to the Orthogonal group case as its known that $\mathrm{Wg}_N^{\mathrm{SP}}(\pi)= \pm \Wg_{2N}^{\mathrm{O}}(\pi)$.
\end{remark}

\subsection{Parings and the Weingarten Graph}

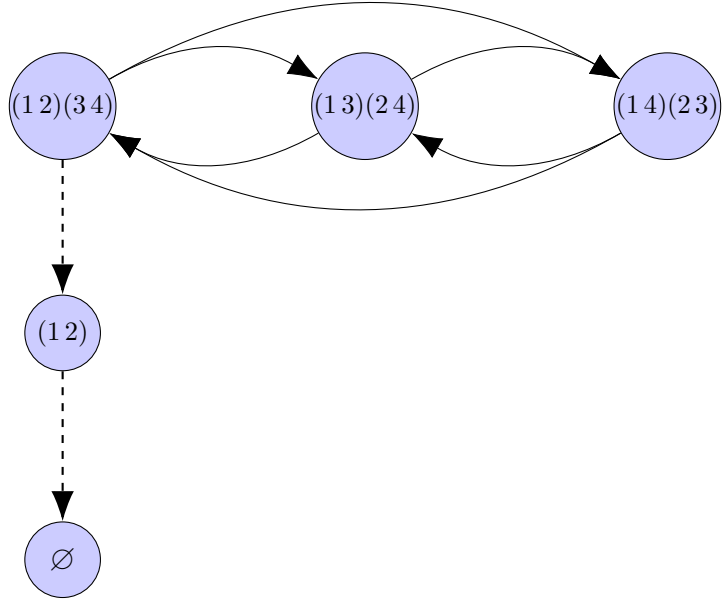
\begin{figure}
    \centering
    \begin{tikzpicture}[
    ->, % Arrow style for directed edges
    >={Latex[length=4mm, width=3mm]}, % Larger arrowheads
    node/.style={circle, draw, fill=blue!20, minimum size=10mm, inner sep=0pt, align=center},
    every edge/.style={draw, thick}
] % Base style for edges
    % Nodes
    \node[node](p1) at (0,6) {$(1 \, 2)(3 \, 4)$};
    \node[node](p2) at (4,6) {$(1\, 3)(2 \, 4)$};
    \node[node](p3) at (8,6) {$(1\, 4)(2\,3)$};
    \node[node](pp1) at (0,3) {$(1\, 2)$};
    \node[node](ppp1) at (0,0) {$\emptyset$};
    \draw (p1) edge[dashed] (pp1);
    \draw (pp1) edge[dashed] (ppp1);
    \draw (p1) to[bend left=30] (p2);
    \draw (p1) to[bend left=30] (p3);
    \draw (p2) to[bend left=30] (p1);
    \draw (p3) to[bend left=30] (p1);
    \draw (p2) to[bend left=30] (p3);
    \draw (p3) to[bend left=30] (p2);
\end{tikzpicture}

\caption{The figure above displays the orthogonal Weingarten graph $\mathcal{W}^{\mathrm{O}}$ restricted to layers 0 through 2, $\bigsqcup_{n=0}^{2} \mathcal{P}_{2}(2n)$.}
\label{fig:Orthogonal Weingarten Graph}
\end{figure}
\begin{definition}
    For a pairing $\pi=\{\{k_1,k_2\},...,\{k_{2n-1},k_{2n}\}\} \in \mathcal{P}_2(2n)$, and a permutation $\sigma \in S_{2n}$, we define the operation, 
    \begin{equs}
        \sigma .\pi := \{(\sigma(k_1),\sigma(k_2)),...,(\sigma(k_{2n-1}),\sigma(k_{2n}))\}
    \end{equs}
\end{definition}

We can define the orthogonal Weingarten graph $\mathcal{W}^{\mathrm{O}}$ to be the directed graph with vertex set $V= \bigsqcup_{n \geq 0} \mathcal{P}_{2}(2n)$. Moreover if $\pi_1, \pi_2 \in \mathcal{P}_{2}(2n)$ and $\pi_2 = (i \hspace{1.5mm} 2n-1). \pi_1$ for some $i<2n-1$, then there is a solid arrow $\pi_1 \to \pi_2$. And if $\pi\in \mathcal{P}_{2}(2n)$ contains the pair $\{2n-1 \hspace{1.5mm} 2n\}$, then there is a dashed arrow going from $\pi$ to $\Tilde{\pi} \in \mathcal{P}_{2}(2n-2)$ which is $\pi$ with the pair $\{2n-1 \hspace{1.5mm} 2n\}$ removed. See Figure \ref{fig:Orthogonal Weingarten Graph} for an image of the orthogonal Weingarten graph

As in the unitary case, for $\pi \in \mathcal{W}^{\mathrm{O}}$ we can define the collection of paths $\mathrm{P}(\pi,l)$ in $\mathcal{W}^{\mathrm{O}}$ going from $\pi \to \emptyset$ using exactly $l$ solid edges. 

We will use the following path expansion proved in \cite{Collins2017} as the definition of $\Wg_N^O (\pi)$ for $n < N$,

\begin{lemma}\label{Orthogonal Weingarten Expansion}
    For $2n \leq N$ and all $\pi \in \mathcal{P}_2(2n)$,
    \begin{equs}
        (-1)^{|\pi|}N^{n+|\pi|}\Wg_N^{\mathrm{O}}(\pi)=\sum_{g \geq 0} |\mathrm{P}(\pi,|\pi|+g)|(-N)^g.
    \end{equs}
    Moreover the series is absolutely convergent. 
\end{lemma}

\begin{remark}
    The absolute convergence of the series in \eqref{Orthogonal Weingarten Expansion} is not mentioned in \cite{Collins2017}, however this can be proven via the trivial bound $|\mathrm{P}(\pi,|\pi|+g)| \leq (2n)^{n+|\pi|+g}$ since there are at most $2n$ options for every edge traversed in a path $\mathrm{p} \in \mathrm{P}(\pi,|\pi|+g)$.
\end{remark}

To each pairing $\pi \in \mathcal{P}_2([2n])$ we can assign a partition as follows.
\begin{definition}
    The partition $\lambda \vdash n$ associated to a pairing $\pi \in \mathcal{P}_2(2n)$ called the \textit{coset type} of $\pi$ can be constructed as follows:
    \begin{enumerate}
        \item Consider a graph with $2n$ vertices labeled $1$ through $2n$.
        \item For each $i \in [n]$, connect vertex $2i-1$ to $2i$ with an edge.
        \item For each pair $\{k_1,k_2\} \in \pi$, connect vertex $k_1$ to $k_2$.
        \item The connected components of this graph clearly have even number of vertices. If the connected components of the graph constructed have sizes $2\lambda_1 \geq ... \geq 2\lambda_l$. Then $\lambda=(\lambda_1,...,\lambda_l) \vdash n$ is the corresponding partition.
    \end{enumerate}
\end{definition}

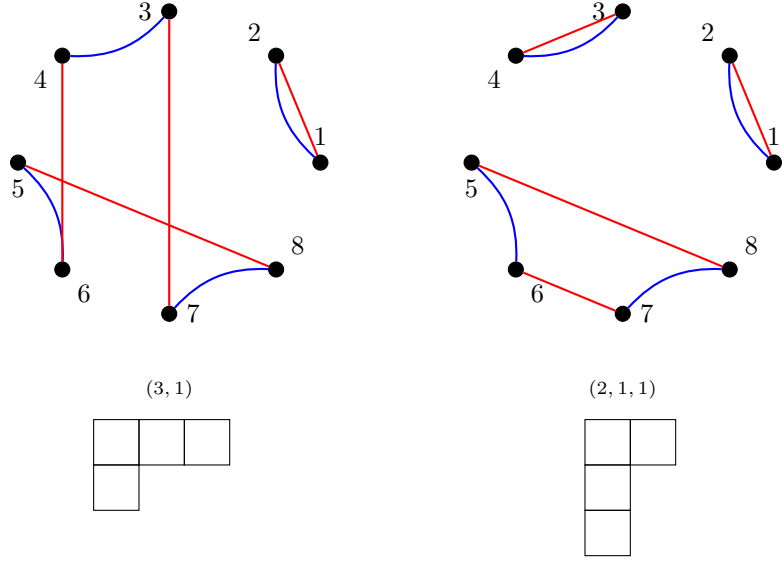
\begin{figure}
    \centering
    \begin{tikzpicture}
    % Left Picture: second to last
    \begin{scope}
        \def\r{2}
        % Store node positions explicitly
        \foreach \i [count=\j from 1] in {0,45,90,135,180,225,270,315} {
            \node[draw, circle, fill=black, inner sep=2pt, label={\i+90:\j}] (N\j) at ({\r*cos(\i)},{\r*sin(\i)}) {};
        }

        % Draw blue arcs between (2i-1) and (2i) using control points for curvature
        \draw[blue, thick] (N1) to[bend left=25] (N2);
        \draw[blue, thick] (N3) to[bend left=25] (N4);
        \draw[blue, thick] (N5) to[bend left=25] (N6);
        \draw[blue, thick] (N7) to[bend left=25] (N8);

        % Draw red edges for specified pairs
        \draw[red, thick] (N1) -- (N2);
        \draw[red, thick] (N3) -- (N7);
        \draw[red, thick] (N4) -- (N6);
        \draw[red, thick] (N5) -- (N8);
    \end{scope}

    % Move the second to last picture to the right
    \node at (5,0) {};  % Adjusting space between the pictures

    % Right Picture: last picture
    \begin{scope}[xshift=6cm] % Shifted to the right
        \def\r{2}
        % Store node positions explicitly
        \foreach \i [count=\j from 1] in {0,45,90,135,180,225,270,315} {
            \node[draw, circle, fill=black, inner sep=2pt, label={\i+90:\j}] (N\j) at ({\r*cos(\i)},{\r*sin(\i)}) {};
        }

        % Draw blue arcs between (2i-1) and (2i) using control points for curvature
        \draw[blue, thick] (N1) to[bend left=25] (N2);
        \draw[blue, thick] (N3) to[bend left=25] (N4);
        \draw[blue, thick] (N5) to[bend left=25] (N6);
        \draw[blue, thick] (N7) to[bend left=25] (N8);

        % Draw new red edges for specified pairs
        \draw[red, thick] (N1) -- (N2);
        \draw[red, thick] (N3) -- (N4);
        \draw[red, thick] (N6) -- (N7);
        \draw[red, thick] (N5) -- (N8);
    \end{scope}

    % Draw Young diagram for the partition (3,1) below the left image
    \node at (0,-3) {\scriptsize \((3,1)\)};
    \draw (-1,-4) rectangle (-0.4,-3.4); % Cell (1,1)
    \draw (-0.4,-4) rectangle (0.2,-3.4); % Cell (1,2)
    \draw (0.2,-4) rectangle (0.8,-3.4); % Cell (1,3)
    \draw (-1,-4.6) rectangle (-0.4,-4);   % Cell (2,1)

    % Draw Young diagram for the partition (2,1,1) below the right image
    \node at (6,-3) {\scriptsize \((2,1,1)\)};
    \draw (5.5,-4) rectangle (6.1,-3.4); % Cell (1,1)
    \draw (6.1,-4) rectangle (6.7,-3.4); % Cell (1,2)
    \draw (5.5,-4.6) rectangle (6.1,-4);   % Cell (2,1)
    \draw (5.5,-5.2) rectangle (6.1,-4.6); % Cell (3,1)

\end{tikzpicture}

\caption{Take $n = 8$, and consider $\pi = (1 \, 2)(3 \, 7)(4 \, 6)(5 \, 8)$, and $\tau = (4 \, 7)$ so that $\tau.\pi=(1 \, 2)(3 \, 4)(6 \, 7)(5 \, 8)$. Then $\mathcal{G}_{\pi}$ is depicted on the left with a Young diagram of the coset type of $\pi$ below. Similarly, $\mathcal{G}_{\tau.\pi}$ is depicted on the right with the Young Diagram of the coset type of $\tau.\pi$ below.}
\label{fig:Pairing Diagram Figure}
\end{figure}

We can interpret the partition defined above diagrammatically as follows. Consider an undirected (multi)graph with vertex set $[2n]$. For each $i \in [n]$, connect $2i-1$ to $2i$ with a blue edge. For each $j \in [2n]$, connect $j$ and $\pi(j)$ with a red edge. The connected components of the resulting graph $\mathcal{G}_{\pi}$ correspond to the coset type of $\pi$. In particular the sizes of the connected components in $\mathcal{G}_{\pi}$ divided by $2$, is precisely the coset type partition of $\pi$. See examples in Figure \ref{fig:Pairing Diagram Figure}.

For a pairing $\pi \in \mathcal{P}_2(2n)$, we can define the norm $|\pi|:=n-\ell(\lambda)$ where $\lambda$ is the partition corresponding to $\pi$ and $\ell(\lambda)$ denotes the number of parts in $\lambda$.

Letting $\mathfrak{e}:=(1 \, 2) (3 \, 4)\cdots (2n-1 \, 2n)$, and interpreting a pairing $\pi \in \mathcal{P}_2(2n)$ as a permutation sending each $i$ to its pair, then the coset type of $\pi$ precisely corresponds to the cycle type of $\pi\mathfrak{e}$. As a consequence we have the following lemma.

\begin{lemma}\label{Orthogonal Transposition Lemma}
    For any $\pi \in \mathcal{P}_2(2n)$ and $\tau = (i \, 2n-1)$ for some $i<2n-1$, then,
    \begin{equs}
        \ell(\tau.\pi)-\ell(\pi) \in \{-1,0,1\}.
    \end{equs}
    Moreover suppose the connected component of $2i-1$ in $\mathcal{G}_{\pi}$ has size $2m$. Then there are exactly $m-1$ choices for $\tau$ such that the coset type of $\tau.\pi$ is the same as $\pi$. And for each $i \in [m-1]$, there is exactly one choice for $\tau$ such that the operation $\pi \to \tau.\pi$ splits the connected component containing $2n-1$ into a component of size $2i$ containing $2n-1$ and a disjoint component of size $2m-2i$.
\end{lemma}
\begin{proof}
    The connected components of $\mathcal{G}_{\pi}$ are just cycles with edges alternating between red and blue color. Let's enumerate the vertices in the cycle containing $2n-1$, by walking around the cycle starting at $2n-1$, and going in the red edge direction to obtain $(2n-1,i_1,...,i_{2m-2},i_{2m-1})$ with $i_{2m-1} = 2n$.

    Case 1: Suppose $\tau = (i \, 2n-1)$ where $i \notin \{i_1,\dots,i_{2m-2}\}$. Then we enumerate the vertices in the cycle containing $i$ in the same way starting at $i$ and walking in the red edge direction, we obtain a sequence $(i, j_1,\dots,j_{m'})$. Then the operation $\pi \to\tau.\pi$ will merge the two cycles $(2n-1,i_1,...,i_{2m-2},i_{2m-1})$ and $(i, j_1,\dots,j_{m'})$ to form the cycle $(2n-1,j_1,j_2,\dots,j_{m'},i,i_{2m-1},i_{2m-2},...,i_1)$. All other cycles will be left untouched. Thus $\ell(\tau.\pi)=\ell(\pi)-1$.

    Case 2: Suppose $\tau = (i_{1} \, 2n-1)$. Then $\tau.\pi=\pi$.
    
    Case 3: Suppose $\tau = (i_{2k+1} \, 2n-1)$ for $k \in [m-2]$. Then the cycle $(2n-1,i_1,...,i_{2m-2},i_{2m-1})$, gets transformed into $(2n-1,i_{2k},i_{2k-1},\dots,i_1,i_{2k+1},i_{2k+2},\dots,i_{2m-1})$ by the operation $\pi \to\tau.\pi$, and so $\ell(\tau.\pi)=\ell(\pi)$, and $\tau.\pi$ and $\pi$ have the same coset type.

    Case 4:  Suppose $\tau = (i_{2k} \, 2n-1)$ for $k \in [m-1]$. Then the cycle $(2n-1,i_1,...,i_{2m-2},i_{2m-1})$, gets split into the two cycles $(2n-1,i_{2k+1},i_{2k+2},\dots,i_{2m-1})$ and $(i_1,i_2,\dots,i_{2k})$ by the operation $\pi \to\tau.\pi$, and so $\ell(\tau.\pi)=\ell(\pi)-1$, and the size of the connected component containing $2n-1$ will now have size $m-2k$.
\end{proof}

As a consequence of the last lemma $\mathrm{P}(\pi,\ell)=\emptyset$ for $\ell< |\pi|$. Moreover, using the previous lemma and the the Catalan number recursion, we obtain the enumeration formula.

\begin{lemma}\label{Orthogonal Walk Enumeration}
    For any pairing $\pi \in \mathcal{P}_2(2n)$ with cycle type $\lambda$,
    \begin{equs}
        |\mathrm{P}(\pi,|\pi|)| = |\mathrm{Moeb}(\lambda)|.
    \end{equs}
\end{lemma}

For convenience we will now reformulate the path expansion Lemma \ref{Orthogonal Weingarten Expansion}. 

\begin{definition}
    First for a path $(\pi_0,...,\pi_{n+\ell}) \in \mathrm{P}(\pi,\ell)$, we will call a solid edge step $\pi_i \to \pi_{i+1}$ a \textit{minor defect} (\textit{major defect}) if $|\pi_{i+1}|=|\pi_i|$ (if $|\pi_{i+1}|=|\pi_i|+1$). Define $\mathrm{P}(\pi,g_1,g_2)$ as the paths in  $\mathrm{P}(\pi,|\pi|+g_1+2g_2)$ with exactly $g_1$ minor defects and $g_2$ major defects. 
\end{definition}

We now state a corollary of Lemmas \ref{Orthogonal Weingarten Expansion} and \ref{Orthogonal Transposition Lemma},

\begin{cor}\label{Refined O Weingarten Expansion}
    For $2n \leq N$ and any pairing $\pi \in \mathcal{P}_2(2n)$,
    \begin{equs}
        (-1)^{|\pi|}N^{n+|\pi|}\mathrm{Wg}_N^{\mathrm{O}}(\pi)= \sum_{g_1,g_2 \geq 0} |\mathrm{P}(\pi,g_1,g_2)| (-N)^{g_1+2g_2}.
    \end{equs}
\end{cor}
\subsection{The Orthogonal Weingarten Process}

We are now in a position to define the orthogonal Weingarten process.

\begin{definition}
    For a partition $\mu \vdash n$, let $\mathcal{P}_2(\mu)$ be the set of pairings $\pi \in\mathcal{P}_2(2n)$ with coset type $\mu$.
\end{definition}

\begin{definition}
    The orthogonal Weingarten process with initial data $\pi\in \mathcal{P}_2(2n)$, $\mathrm{WP}(\pi)$ generates a random path $\mathrm{P}(\pi, |\pi|)$ with the following algorithm.
    \begin{enumerate}
        \item Let $\pi_0 =\pi$, and $n_0=n$. Set $k=0$
        \item If: $\pi_k$ does not contain the pair $\{2n_{k}-1,2n_k\}$, move to the next step. Else if: $n_k>1$ and $\pi_k$ contains the pair $\{2n_{k}-1,2n_k\}$, set $n_{k+1}=n_k-1$, and set $\pi_{k+1} \in \mathcal{P}_2(2n_{k+1})$ to be the pairing $\pi$ with $\{2n_{k}-1,2n_k\}$ removed. Next repeat this step with $k$ now set to $k+1$. Else if: $n_k=1$, terminate the process.
        \item Suppose the cycle in $\mathcal{G}_{\pi}$ containing $2n-1$, is enumerated by a walk starting at $2n-1$, followed by going in the red edge direction as $(2n-1,i_1,i_2,\dots,i_{2m-1})$. For each $j \in [m-1]$, with probability $\frac{\mathrm{Cat(j-1)\mathrm{Cat}(m-j)}}{\mathrm{Cat}(m-1)}$, let $\pi_{k+1}=(i_{2j} \, 2n-1) .\pi_k$, set $n_{k+1}=n_k$, and now set $k$ to $k+1$. Now go back to step 2.
    \end{enumerate}
\end{definition}

\begin{definition}
    For a partition $\mu \vdash n$, let $\mathrm{WP}(\mu)$ be a process generating a random path of pairings $(\pi_0,\pi_1,...,\pi_{2n-\ell(\mu)})$ as follows:
    \begin{enumerate}
        \item Take $\pi \sim \mathrm{Uniform}(\mathcal{P}_2(\mu))$.
        \item Generate a random path $(\pi_0,..., \pi_{2n-\ell(\mu)})$ sampled according to $\mathrm{WP}(\mu)$.
    \end{enumerate}
    Defining $\mu_i := \mathrm{shape}(\pi_i)$, from $\mathrm{WP}(\mu)$, we also obtain a random sequence of partitions, $(\mu_0,\mu_1,...,\mu_{2n-\ell(\lambda)})$.
\end{definition}

$\mathrm{WP}(\pi)$ and $\mathrm{WP}(\mu)$ come with the canonical filtration $(\mathcal{F}_k)_{k \leq n+|\pi|}$ where $\mathcal{F}_k$ is the sigma algebra generated by $\{\pi_j\}_{j \leq k}$. Moreover, $\mathrm{WP}(\mu)$ has another natural filtration, $(\Sigma_k)_{k \leq 2n-\ell(\lambda)}$, where $\Sigma_k$ is the $\sigma$-algebra generated by $\{\mu_j\}_{j \leq k}$.

\begin{lemma}\label{Orthogonal Process is uniform over paths}
    For any permutation $\pi \in \mathcal{P}_2(2n)$,
    \begin{equs}
        \mathrm{WP}(\pi) \sim \mathrm{Uniform}(\mathrm{P}(\pi, |\pi|)),
    \end{equs}
    and for any $\mu \vdash n$,
    \begin{equs}
        \mathrm{WP}(\mu) \sim \mathrm{Uniform}\bigg( \bigcup_{\pi \in \mathcal{P}_2(\mu)} \mathrm{P}(\pi, |\pi|)\bigg).
    \end{equs}
\end{lemma}

\begin{proof}

    The proof proceeds by induction on $|\pi|$ as in the case of Lemma \ref{Unitary Process is uniform over paths}. The base case is trivial

    For the inductive step, suppose first that $\pi$ does not contain the pair $(2n-1,2n)$. Then let $(2n-1,i_1,i_2,\dots,i_{2m-1})$ be the cycle containing $2n-1$ in $\mathcal{G}_{\pi}$. Then by Lemma \ref{Orthogonal Transposition Lemma} the first step in $\mathrm{P}(\pi,|\pi|)$ must take $\pi \to (i_{k} \, 2n-1).\pi$ for some $k$. Thus
    \begin{equs}
        \mathrm{P}(\pi,\ell+1,0) \simeq \bigcup_{k=1}^{m-1}\mathrm{P}((i_{2k} \, 2n-1).\pi,|\pi|).
    \end{equs}
    Moreover, using the enumeration formula in Lemma \ref{Orthogonal Walk Enumeration}, sampling $p \in \mathrm{P}(\pi,|\pi|)$ uniformly, with probability
    \begin{equs}
        \frac{|\mathrm{P}((i_{2k} \, 2n-1).\pi, \ell+1,0)|}{|\mathrm{P}(\pi,\ell+1,0)|}=\frac{\mathrm{Cat(k-1)\mathrm{Cat}(m-k)}}{\mathrm{Cat}(m-1)}
    \end{equs}
    the first step of the path is $\pi \to (i_{2k} \, 2n-1).\pi$, followed by a uniformly random path in $\mathrm{P}((i_{2k} \, 2n-1).\pi, |\pi|-1)$. But the Weingarten process obeys this recursion by construction.
\end{proof}

At this point we remark that the partition $\mu_j$ at the $j$th step of the orthogonal $\mathrm{WP}^{\mathrm{O}}(\mu)$ has exactly the same law as the unitary Weingarten process $\mathrm{WP}^{\mathrm{U}}(\mu)$ since the transition probabilities defining the processes are exactly the same!

As in Subsection \ref{Reduction to Expectation Section}, we will use $L_j$ to denote the length of the longest part of $\mu_j$, the partition at the $j$th step in Weingarten process $\mathrm{WP}^{\mathrm{O}}(\mu)$ at the. We once again remark that this $L_j$ has the same distributions as in the unitary case.

\subsection{Proof of Theorem \ref{Main Bound Orthogonal}}
For a path $\mathrm{p} \in \mathrm{P}(\pi,g_1,g_2)$ we can define $j:= j(p)$ such that the $(j+1)$th solid arrow crossed in the path is the first transition $\pi_{j+n-r} \to \pi_{j+n-r+1}$ such that $|\pi_{j+n-r}|\leq|\pi_{j+n-r+1}|$, and we define
\begin{equs}
    \mathrm{P}_j(\pi,g_1,g_2):=\{p \in \mathrm{P}(\pi,g_1,g_2): j(p)=j\}.
\end{equs}
We clearly have
\begin{equs}\label{First Backtack Decomp Orthogonal}
    \mathrm{P}(\pi,g_1,g_2)= \bigcup_{j=0}^{|\sigma|} (\mathrm{P}_j(\pi,g_1-1,g_2) \cup \mathrm{P}_j(\pi,g_1,g_2-1)). 
\end{equs}

Moreover, every path can be decomposed into the path up to this step, and the remaining part: i.e. for $|\pi|=k$,
\begin{equs}\label{Orthogonal First Backtack, level Decomp}
    \mathrm{P}_j(\pi,g_1,g_2) = \bigcup_{r=2}^{k} \bigcup_{\eta \in \mathcal{P}_2(2r)} \bigcup_{\substack{\eta' \in \mathcal{P}_2(2r), |\eta'|=|\eta|,\eta \to \eta' \\ \eta' \in \mathcal{P}_2(2r), |\eta'|=|\eta|+1,\eta \to \eta'}} \mathrm{P}_j(\pi,\eta,\eta',g_1,g_2) ,
\end{equs}
where $\mathrm{P}_j(\pi,\eta,\eta',g_1,g_2)$ denotes the set of paths $\mathrm{P}_j(\pi,g_1,g_2)$, where the first ste $\pi_{j+n-r} \to \pi_{j+n-r+1}$ with $|\pi_{j+n-r}|\leq |\pi_{j+n-r+1}|$ occurs at $\pi_{j+n-r} =\eta$ and $\pi_{j+k-r+1}=\eta'$.

\noindent If we let
\begin{equs}
    \Tilde{\mathrm{P}}_j(\pi,\eta),
\end{equs}
denotes the set of all "partial paths" from $\pi \to \eta$, i.e. paths $\pi= \pi_0, \pi_1,...,\pi_{j+n-r}=\eta$ with no increasing norm steps. Then for $|\eta|=|\eta'|$ we have the factorization,

\begin{equs}\label{Orthogonal Partial path decomposition 1}
    \mathrm{P}_j(\pi,\eta,\eta',g_1,g_2) \simeq \Tilde{\mathrm{P}}_j(\pi,\eta) \times \mathrm{P}(\eta', g_1-1,g_2).
\end{equs}
And for $|\eta|=|\eta'|$,
\begin{equs}\label{Orthogonal Partial path decomposition eq 2}
    \mathrm{P}_j(\pi,\eta,\eta',g_1,g_2) \simeq \Tilde{\mathrm{P}}_j(\pi,\eta) \times \mathrm{P}(\eta', g_1,g_2-1).
\end{equs}

The next lemma is a corollary of this path decomposition.

\begin{lemma}\label{Orthogonal Expectation Reduction}
    Suppose $\pi$ has coset type $\mu$, then we have,
    \begin{equs}\label{Minor Defect Path Bound}
        \frac{|\mathrm{P}(\pi,1,0)|}{|\mathrm{P}(\pi,0,0)|} \leq \sum_{j=1}^{n+|\pi|} \E_{\mathrm{WP}(\mu)}[L_j],
    \end{equs}
    and
    \begin{equs}\label{Major Defect Path Bound}
        \frac{|\mathrm{P}(\pi,0,1)|}{|\mathrm{P}(\pi,0,0)|} \leq 6\sqrt{8}n\sum_{j=1}^{n+|\pi|} \E_{\mathrm{WP}(\mu)}[L_j^{3/2}].
    \end{equs}
\end{lemma}
\begin{proof}
    We only prove \eqref{Minor Defect Path Bound}, as \eqref{Major Defect Path Bound} is proven identically to Lemma \ref{Single Defect Path Bound}. Combining \eqref{Orthogonal First Backtack, level Decomp} and \eqref{Orthogonal Partial path decomposition 1},
    \begin{equs}
        |\mathrm{P}(\pi,1,0)| &= \sum_{j=0}^{|\pi|} |\mathrm{P}_j(\pi, 1,0)|\\
        &= \sum_{j=0}^{|\pi|} \sum_{r=2}^{n} \sum_{\eta \in \mathcal{P}_2(2r)} \sum_{\eta' \in \mathcal{P}_2(2r) \backslash \{\eta\}, |\eta'|=|\eta|,\eta \to \eta'} |\mathrm{P}_j(\pi,\eta,\eta',1,0)| \\
        &= \sum_{j=0}^{|\pi|} \sum_{r=2}^{n} \sum_{\eta \in \mathcal{P}_2(2r)}\sum_{\eta' \in \mathcal{P}_2(2r) \backslash \{\eta\}, |\eta'|=|\eta|,\eta \to \eta'}|\Tilde{\mathrm{P}}_j(\pi,\eta)| |\mathrm{P}(\eta',0,0)|.  
    \end{equs}

Now, recall that if $\eta \to \eta'$ with $|\eta|=|\eta'|$, then by Lemma \ref{Orthogonal Transposition Lemma}, $\eta$ and $\eta'$ have the same coset type, and hence, $\mathrm{P}(\eta',0,0)=\mathrm{P}(\eta,0,0)=\mathrm{P}(\eta,|\eta|)$. As a consequence, letting $\ell(\pi)$ denoting the size of the largest part of the coset type of $\pi$,
\begin{equs}
        &|\mathrm{P}(\pi,1,0)|\\
        &=\sum_{j=0}^{|\pi|} \sum_{r=2}^{n} \sum_{\eta \in \mathcal{P}_2(2r)}\sum_{\eta' \in \mathcal{P}_2(2r) \backslash \{\eta\}, |\eta'|=|\eta|,\eta \to \eta'}|\Tilde{\mathrm{P}}_j(\pi,\eta)| |\mathrm{P}(\eta,|\eta|)|\\
        &\leq \sum_{j=0}^{|\pi|} \sum_{r=2}^{n} \sum_{\eta \in \mathcal{P}_2(2r)}|\Tilde{\mathrm{P}}_j(\pi,\eta)| |\mathrm{P}(\eta,|\eta|)| \ell(\pi_{n-r+j})\\
        &\leq |\mathrm{P}(\pi,|\pi|)|\sum_{j=0}^{n+|\pi|}  \E_{\mathrm{WP}(\pi)}[L_{j}].
\end{equs}

\end{proof}

From the last Lemma, together with Proposition \ref{prop:Main expectation bound}, we are ready to show the following. 

\begin{lemma}\label{Orthogonal Genus 1 Path Count Bound}
    For any $\pi \in \mc{P}_2(2n)$, we have ,
    \begin{equs}\label{eq:Orth-Path-Defect-1}
        \frac{|\mathrm{P}(\pi,1,0)|}{|\mathrm{P}(\pi,0,0)|} \leq 10^6n^{3/2},
    \end{equs}
    and
    \begin{equs}\label{eq:Orth-Path-Defect-2}
        \frac{|\mathrm{P}(\pi,0,1)|}{|\mathrm{P}(\pi,0,0)|} \leq 6\sqrt{8} \times 10^6n^{3}.
    \end{equs}
\end{lemma}

\begin{proof}
    Since $L_j$ has the same distribution for the orthogonal Weingarten process has the same distribution as for the unitary Weingarten process, \eqref{eq:Orth-Path-Defect-2} is an immediate corollary of Proposition \ref{prop:Main expectation bound} together with Lemma \ref{Orthogonal Expectation Reduction}.
    
For \eqref{eq:Orth-Path-Defect-1}, we argue similarly to the proof of Proposition \ref{prop:Main expectation bound}, applying strong induction on $n$ to prove the slightly stronger claim that $\frac{|\mathrm{P}(\pi,1,0)|}{|\mathrm{P}(\pi,0,0)|} \leq 10^6n^{}L_0^{1/2}$ by using the unitary Weingarten process. 

For the base case when $n \leq 6$ or when $L_0 < 6$, the claim is obvious and there is nothing to prove.

For the inductive step, we split the sum $\sum_{j=0}^{n+|\sigma|} \E[L_j]$ into two parts at the index $T$, bounding the $L_j$ until $T$ by $L_0$, and applying the bound of Lemma \ref{lm:time-to-halve} to the first part, while applying the Markov property, \eqref{eq:general-pivotal-markov}, and the inductive hypothesis to the second part as follows,
    \begin{equs}
        ~&\frac{|\mathrm{P}(\pi,1,0)|}{|\mathrm{P}(\pi,0,0)|}\\
        &\leq \sum_{j=0}^{n+|\sigma|} \E[L_j] \\
        &=\E\bigg[\sum_{j=0}^{T-1} L_j\bigg] +\E\bigg[\sum_{j=T}^{n+|\sigma|} L_j\bigg]\\
        &\leq \E[T]L_0+\E_{\mrm{WP}(\lambda)}\bigg[\E_{\mrm{WP}(\lambda_T)}\bigg[\sum_{j=0}^{n+|\sigma_T|} L_j\bigg]\bigg]\\
        &\leq 10^4nL_0^{1/2}\bigg(1+\log\bigg(\frac{3n}{2L_0}\bigg)\bigg)^2+ \E_{\mrm{WP}(\lambda)}[10^6 nL_T^{1/2}]\\
        &\leq 10^4nL_0^{1/2} \times 4 (3n/2L_0)^{1/2}+(2/3)^{1/2}10^6 n L_0^{1/2}\\
        &\leq 10^6 nL_0^{1/2}
    \end{equs}
    where we used the elementary inequalities $1+\log(x) \leq 4x^{1/4}$ for $x \geq 1$, and $4\times (3/2)^{1/2}+100(2/3)^{1/2}\leq 100$.
    %where we used the fact that $L_T \leq 2L_0/3$ by definition. Lastly the proposition follows from the inequality $n x\log(n/x) \leq n^2$ for $0<x \leq n$ which can be verified by elementary calculus.
\end{proof}

Now via induction we have the following general path counting bound.
\begin{lemma}\label{Orthogonal Main Path Count Bound}
    For any $\pi \in \mc{P}_2(2n)$,
    \begin{equs}
        \frac{|\mathrm{P}(\pi,g_1,g_2)|}{|\mathrm{P}(\pi,0,0)|} \leq (10^6n)^{3g_1/2+3g_2}.
    \end{equs}
\end{lemma}

\begin{proof}
    The case $g_1+g_2 \leq 1$ is the statement of the previous lemma. The general case now follows by induction: 

    \begin{equs}
        |\mathrm{P}(\pi,g_1,g_2)| &= \sum_{j=0}^{|\pi|} |\mathrm{P}_j(\pi, g_1,g_2)|\\
        &= \sum_{j=0}^{|\pi|} \sum_{r=2}^{n} \sum_{\eta \in \mathcal{P}_2(2r)} \sum_{\eta' \in \mathcal{P}_2(2r) \backslash \{\eta\}, |\eta'|=|\eta|,\eta \to \eta'} |\mathrm{P}_j(\pi,\eta,\eta',g_1,g_2)| \\
        &= \sum_{j=0}^{|\pi|} \sum_{r=2}^{n} \sum_{\eta \in \mathcal{P}_2(2r)}\sum_{\eta' \in \mathcal{P}_2(2r) \backslash \{\eta\}, |\eta'|=|\eta|,\eta \to \eta'}|\Tilde{\mathrm{P}}_j(\pi,\eta)| |\mathrm{P}(\eta',g_1-1,g_2)|\\
        &+\sum_{j=0}^{|\pi|} \sum_{r=2}^{n} \sum_{\eta \in \mathcal{P}_2(2r)}\sum_{\eta' \in \mathcal{P}_2(2r) \backslash \{\eta\}, |\eta'|=|\eta|+1,\eta \to \eta'}|\Tilde{\mathrm{P}}_j(\pi,\eta)| |\mathrm{P}(\eta',g_1,g_2-1)|\\
        &\leq \sum_{j=0}^{|\pi|} \sum_{r=2}^{n} \sum_{\eta \in \mathcal{P}_2(2r)}\sum_{\eta' \in \mathcal{P}_2(2r) \backslash \{\eta\}, |\eta'|=|\eta|,\eta \to \eta'}(Cn)^{3(g_1-1)/2+3g_2}|\Tilde{\mathrm{P}}_j(\pi,\eta)| |\mathrm{P}(\eta',0,0)|\\
        &+\sum_{j=0}^{|\pi|} \sum_{r=2}^{n} \sum_{\eta \in \mathcal{P}_2(2r)}\sum_{\eta' \in \mathcal{P}_2(2r) \backslash \{\eta\}, |\eta'|=|\eta|+1,\eta \to \eta'}(Cn)^{3g_1/2+3(g_2-1)}|\Tilde{\mathrm{P}}_j(\pi,\eta)| |\mathrm{P}(\eta',0,0)|\\
        &= (Cn)^{3(g_1-1)/2+3g_2} |\mathrm{P}(\pi,1,0)|+(Cn)^{3g_1/2+3(g_2-1)}|\mathrm{P}(\pi,0,1)|\\
        &\leq (Cn)^{3g_1/2+3g_2}|\mathrm{P}(\pi,0,0)|.
    \end{equs}
\end{proof}

Finally Theorem \ref{Main Bound Orthogonal} follows from the prior path counting bound.

\begin{proof}[Proof of Theorem \ref{Main Bound Orthogonal}]
    Substituting in Lemma \ref{Orthogonal Main Path Count Bound} into the path expansion Corollary \ref{Refined O Weingarten Expansion}, setting $C=10^6$ we have,
\begin{equs}
    N^{n+|\pi|}|\Wg_N^{\mathrm{O}}(\pi)| &\leq \sum_{g_1,g_2 \geq 0} |\mathrm{P}(\pi,g_1,g_2)|N^{-g_1-2g_2} \\
    &\leq \sum_{g_1,g_2 \geq 0} (Cn^{3/2}/N)^{g_1+2g_2} |\mathrm{Moeb}(\pi)|\\
    & \leq \sum_{g \geq 0} (g+1)(Cn^{3/2}/N)^{g} |\mathrm{Moeb}(\pi)|\\
    &\leq \frac{|\mathrm{Moeb}(\pi)|}{(1-Cn^{3/2}/N)^2},
\end{equs}
where we used the fact that the equation $g_1+2g_2=g$ has $\leq g+1$, and that $\sum_{g \geq 0}(g+1)x^g = (1-x)^{-2}$ for $|x|<1$.
\end{proof}

\vspace{10mm}
\noindent
\textbf{\Large Declarations}

\vspace{3mm}
\noindent
\textbf{Conflict of Interest} The author declares that there is no conflict of interest.

\vspace{3mm}
\noindent
\textbf{Data Availability} We do not analyse or generate any datasets, because our work
proceeds within a theoretical and mathematical approach. One
can obtain the relevant materials from the references below.

\bibliographystyle{alpha}
\bibliography{references}
\Addresses
\end{document}